# ISS IN DIFFERENT NORMS FOR 1-D PARABOLIC PDES WITH BOUNDARY DISTURBANCES

**Iasson Karafyllis[*] and Miroslav Krstic[**]**

[*]Dept. of Mathematics, National Technical University of Athens,
Zografou Campus, 15780, Athens, Greece, email: iasonkar@central.ntua.gr

[**]Dept. of Mechanical and Aerospace Eng., University of California, San Diego, La Jolla, CA 92093-0411, U.S.A., email: krstic@ucsd.edu

**Abstract**
For 1-D parabolic PDEs with disturbances at both boundaries and distributed disturbances we provide ISS estimates in various norms. Due to the lack of an ISS Lyapunov functional for boundary disturbances, the proof methodology uses (i) an eigenfunction expansion of the solution, and (ii) a finite-difference scheme. The ISS estimate for the sup-norm leads to a refinement of the well-known maximum principle for the heat equation. Finally, the obtained results are applied to quasi-static thermoelasticity models that involve nonlocal boundary conditions. Small-gain conditions that guarantee the global exponential stability of the zero solution for such models are derived.

## 1. Introduction

The extension of the Input-to-State Stability (ISS) property (first developed by E. D. Sontag in [27]) to systems which are described by Partial Differential Equations (PDEs) requires novel mathematical tools and novel approaches. Extensions of the ISS property to PDE systems have been pursued in many recent articles (see for example [1,2,4,5,6,,13,14,17,18,19,20,21,24]). In particular, for PDE systems there are two qualitatively distinct locations where a disturbance can appear: the domain (a distributed disturbance appearing in the PDE) and the boundary (a disturbance that appears in the boundary conditions). Most of the existing results in the literature deal with distributed disturbances.

   Boundary disturbances present a major challenge, because the transformation of the boundary disturbance to a distributed disturbance leads to the ISS property with respect to the boundary disturbance and some of its time derivatives (see for example [1]). This is explained by the use of unbounded operators for the expression of the effect of the boundary disturbance (see the relevant discussion in [17] for inputs in infinite-dimensional systems that are expressed by means of unbounded linear operators). Moreover, although the construction of Lyapunov functionals for PDEs has progressed significantly during the last years (see for example [2,12,13,16,18,21,24]), none of the proposed Lyapunov functionals can be used for the derivation of the ISS property w.r.t. boundary disturbances. Indeed, there is no guarantee that an ISS Lyapunov functional exists for PDEs with boundary disturbances, since all converse Lyapunov results require the so-called "concatenation property" for the allowed disturbances (see [12,21]). Unfortunately, the "concatenation property" does not hold for boundary disturbances. This complexity was noted in [14], where the state space was time-varying and closely related to the boundary disturbance itself.

   The recent article [14] suggested a methodology for obtaining ISS estimates w.r.t. boundary disturbances for 1-D parabolic PDEs. The main idea in [14] was the eigenfunction expansion of the solution and the avoidance of requiring the state to take values in a constant state space: the evolution of the state was regarded to take place in a parameterized convex set. The obtained ISS estimates were expressed in weighted $L^2$ norms.

This paper focuses on 1-D parabolic PDEs with disturbances acting on both boundary sides and distributed disturbances. The purpose is the establishment of the ISS property in various norms. While an ISS estimate in the $L^2$ norm is derived (Theorem 2.3) by performing an analysis similar to that proposed in [14], this is not the case with different function norms. To this end, we develop here a finite-difference scheme, which allows the derivation of ISS estimates expressed in weighted $L^\infty$ and $L^1$ norms (Theorem 2.2 and Theorem 2.4). Finite-difference schemes have been used with success in the past for various PDEs (see for example [11]). The obtained estimates of the gains of the boundary disturbances are exact, in the sense that constant disturbances lead to equilibrium profiles for which the norms are exactly equal to the estimated gains.

Therefore, it can be claimed that the present paper develops a novel mathematical tool that can be used for the derivation of the ISS property for boundary disturbances: the utilization of finite-difference schemes. Another important feature of the present work is the fact that the ISS estimate in $L^\infty$ norm for the heat equation (Corollary 2.5) leads to a refinement of the well-known maximum principle for the heat equation (see for example [11] on page 215).

The obtained ISS estimates can be exploited, by means of small-gain arguments, for the stability analysis of parabolic PDEs with nonlocal boundary conditions. This is shown in Section 3 of our paper, where we study a 1-D quasi-static thermoelasticity model. The model involves nonlocal boundary conditions and was proposed and first studied in [7,8]. Stability conditions for 1-D quasi-static thermoelasticity models have been given in [7,8] and in [9,10,15] (for numerical stability of discretization schemes). We exploit the ISS estimates for 1-D parabolic PDEs and provide small-gain conditions for the global exponential stability of the zero solution in various norms (Theorem 3.1, Theorem 3.2 and Theorem 3.3). The obtained stability conditions are weaker than the conditions proposed in the literature.

The structure of the paper is as follows. Section 2 is devoted to the presentation of the problem and the statement of the main results which allow the derivation of the ISS estimates in various norms (Theorem 2.2, Theorem 2.3 and Theorem 2.4). Specific results are provided for the heat equation with Dirichlet boundary conditions (Corollary 2.5). The main results on the global exponential stability in various norms of the zero solution for quasi-static thermoelasticity models are given in Section 3 (Theorem 3.1, Theorem 3.2 and Theorem 3.3). Section 4 of the present work contains the proofs of the main results. The concluding remarks of the paper are provided in Section 5. Finally, the Appendix contains the proofs of some auxiliary technical results, which were used in Section 4.

**Notation.** Throughout this paper, we adopt the following notation.

∗ Let $U \subseteq \Re^n$ be a set with non-empty interior and let $\Omega \subseteq \Re$ be a set. By $C^0(U;\Omega)$, we denote the class of continuous mappings on $U$, which take values in $\Omega$. By $C^k(U;\Omega)$, where $k \geq 1$, we denote the class of continuous functions on $U$, which have continuous derivatives of order $k$ on $U$ and take values in $\Omega$. When $\Omega \subseteq \Re$ is not explicitly given, i.e., when we write $C^k(U)$, we mean that $\Omega = \Re$. For $x \in C^1([0,1])$, $x'(z)$ denotes the derivative with respect to $z \in [0,1]$.

∗ $\Re_+$ denotes the set of non-negative real numbers. For a real number $x \in \Re$, $[x]$ denotes the integer part of $x$, i.e., the greatest integer which is less or equal to $x$.

∗ Let $x \in C^0(\Re_+ \times [0,1])$ be given. We use the notation $x[t]$ to denote the profile at certain $t \geq 0$, i.e., $(x[t])(z) = x(t,z)$ for all $z \in [0,1]$. $H^2(0,1)$ denotes the Sobolev space of continuously differentiable functions on $[0,1]$ with measurable, square integrable second derivative. $L^2_r([0,1])$ denotes the equivalence class of measurable functions $f:[0,1] \to \Re$ for which $\left( \int_0^1 r(z)|f(z)|^2 dz \right)^{1/2} < +\infty$.

## 2. ISS in Various Norms

Let $p \in C^1([0,1];(0,+\infty))$, $r \in C^0([0,1];(0,+\infty))$, $q \in C^0([0,1])$ be given functions and consider the Sturm-Liouville operator $A: D \to L_r^2([0,1])$, defined by

$$(Af)(z) = -\frac{1}{r(z)}\left(p(z)f'(z)\right)' + \frac{q(z)}{r(z)}f(z),$$

for all $f \in D$ and $z \in (0,1)$ \hfill (2.1)

where $D$ is the set of all functions $f \in H^2(0,1)$ for which

$$g_0 f(0) + v_0 f'(0) = g_1 f(1) + v_1 f'(1) = 0 \tag{2.2}$$

where $g_0, v_0, g_1, v_1$ are real constants satisfying the following assumptions:

**(H1)** *Either $v_0 > 0$ or $v_0 = 0$ and $g_0 < 0$.*

**(H2)** *Either $v_1 > 0$ or $v_1 = 0$ and $g_1 > 0$.*

**FACT** (see Chapter 11 in [3] and pages 498-505 in [22]): All eigenvalues of the Sturm-Liouville operator $A: D \to L_r^2([0,1])$, defined by (2.1), (2.2) under (H1), (H2) are real. They form an infinite, increasing sequence $\lambda_1 < \lambda_2 < \ldots < \lambda_n < \ldots$ with $\lim_{n \to \infty}(\lambda_n) = +\infty$. To each eigenvalue $\lambda_n \in \Re$ ($n = 1,2,\ldots$) corresponds exactly one eigenfunction $\phi_n \in D \cap C^2([0,1])$ that satisfies $A\phi_n = \lambda_n \phi_n$. The eigenfunctions form an orthonormal basis of $L_r^2([0,1])$.

In this paper (as in [14]), we also make the following assumption for the Sturm-Liouville operator $A: D \to L_r^2([0,1])$ defined by (2.1), (2.2).

**(H3):** *The Sturm-Liouville operator $A: D \to L_r^2([0,1])$ defined by (2.1), (2.2), where $g_0, v_0, g_1, v_1$ are real constants satisfying (H1), (H2), satisfies*

$$\lambda_1 > 0 \tag{2.3}$$

*and*

$$\sum_{n=1}^{\infty} \lambda_n^{-1} \max_{0 \le z \le 1}\left(|\phi_n(z)|\right) < +\infty \tag{2.4}$$

Consider the parameterized convex set

$$X(\mu_0, \mu_1) := \left\{ x \in C^2([0,1]) : g_0 x(0) + v_0 x'(0) = \mu_0,\ g_1 x(1) + v_1 x'(1) = \mu_1 \right\} \tag{2.5}$$

with parameters $\mu_0, \mu_1 \in \Re$. Given $d_0, d_1 \in C^0(\Re_+)$, $u \in C^1(\Re_+ \times [0,1])$ and $x_0 \in X(d_0(0), d_1(0))$, we study the solution $x \in C^0(\Re_+ \times [0,1]; \Re) \cap C^1((0,+\infty) \times [0,1]; \Re)$ for which $x[t] \in X(d_0(t), d_1(t))$ for all $t \ge 0$, $x(0, z) = x_0(z)$ for all $z \in [0,1]$ and

$$\frac{\partial x}{\partial t}(t, z) = \frac{p(z)}{r(z)} \frac{\partial^2 x}{\partial z^2}(t, z) + \frac{p'(z)}{r(z)} \frac{\partial x}{\partial z}(t, z) - \frac{q(z)}{r(z)} x(t, z) + u(t, z), \text{ for all } (t, z) \in (0, +\infty) \times (0,1) \tag{2.6}$$

In other words, we consider the solution of the evolution equation (2.6) that satisfies for all $t \ge 0$ the boundary conditions

$$g_0 x(t,0) + v_0 \frac{\partial x}{\partial z}(t,0) = d_0(t)$$
$$g_1 x(t,1) + v_1 \frac{\partial x}{\partial z}(t,1) = d_1(t)$$
(2.7)

The inputs $d_0, d_1 \in C^0(\Re_+)$ are boundary disturbances and appear only at the boundary condition, while the input $u \in C^1(\Re_+ \times [0,1])$ is a distributed input that appears in the differential equation. Under sufficient regularity properties for the disturbances $d_0, d_1 \in C^0(\Re_+)$ we are in a position to guarantee that a classical solution exists for the evolution equation (2.6) with (2.7). This is guaranteed by the following result.

**Theorem 2.1:** *Suppose that assumptions (H1), (H2), (H3) hold. Then for every $d_0, d_1 \in C^2(\Re_+)$, $u \in C^1(\Re_+ \times [0,1])$ and $x_0 \in X(d_0(0), d_1(0))$, the evolution equation (2.6) with (2.7) and initial condition $x_0 \in X(d_0(0), d_1(0))$ has a unique solution $x \in C^0(\Re_+ \times [0,1]) \cap C^1((0,+\infty) \times [0,1])$ for which $x[t] \in X(d_0(t), d_1(t))$ for all $t \geq 0$ and $x(0,z) = x_0(z)$ for all $z \in [0,1]$.*

Theorem 2.1 guarantees that for any given $x_0 \in C^2([0,1])$, $u \in C^1(\Re_+ \times [0,1])$, there exists a non-empty set $\Phi(x_0, u) \subseteq C^1(\Re_+; \Re^2)$ such that the following implication holds:

"If $(d_0, d_1) \in \Phi(x_0, u)$ then the evolution equation (2.6) with (2.7) and initial condition $x_0 \in C^2([0,1])$ has a unique solution $x \in C^0(\Re_+ \times [0,1]) \cap C^1((0,+\infty) \times [0,1])$ for which $x[t] \in C^2([0,1])$ for all $t \geq 0$ and $x(0,z) = x_0(z)$ for all $z \in [0,1]$."

Indeed, Theorem 2.1 guarantees that

$$\left\{ (d_0, d_1) \in C^2(\Re_+; \Re^2) : d_0(0) = g_0 x_0(0) + v_0 x_0'(0), d_1(0) = g_1 x_0(1) + v_1 x_0'(1) \right\} \subseteq \Phi(x_0, u)$$

The following result is the first main result of the section. It provides an ISS estimate of the solution of (2.6), (2.7) in a weighted $L^\infty$ norm under the following assumption.

**(H4)** *There exists a function $\eta \in C^2([0,1]; (0,+\infty))$ and a constant $\sigma > 0$ such that $p(z)\eta''(z) + p'(z)\eta'(z) - q(z)\eta(z) = -\sigma r(z)\eta(z)$ for all $z \in [0,1]$. Moreover, the inequalities $g_0 \eta(0) + v_0 \eta'(0) < 0$ and $g_1 \eta(1) + v_1 \eta'(1) > 0$ hold.*

**Theorem 2.2:** *Suppose that assumptions (H1), (H2), (H3), (H4) hold. Then for every $x_0 \in C^2([0,1])$, $u \in C^1(\Re_+ \times [0,1])$ and $(d_0, d_1) \in \Phi(x_0, u)$, the unique solution $x \in C^0(\Re_+ \times [0,1]) \cap C^1((0,+\infty) \times [0,1])$ of the evolution equation (2.6) with (2.7) and initial condition $x_0 \in C^2([0,1])$ satisfies the following estimate for all $t \geq 0$:*

$$\|x[t]\|_{\infty,\eta} \leq \max\left( \exp(-\sigma t) \|x[0]\|_{\infty,\eta}, \frac{\max_{0 \leq s \leq t}(|d_0(s)|)}{|g_0 \eta(0) + v_0 \eta'(0)|}, \frac{\max_{0 \leq s \leq t}(|d_1(s)|)}{g_1 \eta(1) + v_1 \eta'(1)} \right) + \sigma^{-1} \max_{0 \leq s \leq t}\left( \max_{0 \leq z \leq 1}\left( \frac{u(s,z)}{\eta(z)} \right) \right) \quad (2.8)$$

*where*

$$\|x[t]\|_{\infty,\eta} := \max_{0 \leq z \leq 1}\left( \frac{|x(t,z)|}{\eta(z)} \right). \quad (2.9)$$

The second main result of the present work provides ISS estimates of the solution of (2.6), (2.7) in a weighted $L^2$ norm.

**Theorem 2.3:** *Suppose that assumptions (H1), (H2), (H3) hold. Then for every $x_0 \in C^2([0,1])$, $u \in C^1(\Re_+ \times [0,1])$ and $(d_0, d_1) \in \Phi(x_0, u)$, the unique solution $x \in C^0(\Re_+ \times [0,1]) \cap C^1((0,+\infty) \times [0,1])$ of the evolution equation (2.6) with (2.7) and initial condition $x_0 \in C^2([0,1])$ satisfies the following estimate for all $\varepsilon, \omega > 0$ and $t \geq 0$:*

$$\|x[t]\|_{2,r} \leq \sqrt{\frac{\exp(-\lambda_1 t)}{2 - \exp(-\lambda_1 t)}} \|x[0]\|_{2,r} + C_0 \sqrt{(1+\varepsilon^{-1})(1+\omega)} \max_{0 \leq s \leq t}(|d_0(s)|)$$
$$+ C_1 \sqrt{(1+\varepsilon)(1+\omega)} \max_{0 \leq s \leq t}(|d_1(s)|) + \tilde{C}\sqrt{1+\omega^{-1}} \max_{0 \leq s \leq t}\left(\max_{0 \leq z \leq 1}(|u(s,z)|)\right) \quad (2.10)$$

*where*

$$C_0 := \frac{p(0)}{g_0^2 + v_0^2} \sqrt{\sum_{n=1}^{\infty} \frac{1}{\lambda_n^2} \left| g_0 \frac{d\phi_n}{dz}(0) - v_0 \phi_n(0) \right|^2} = \frac{1}{\sqrt{g_0^2 + v_0^2}} \|\tilde{x}\|_{2,r}, \quad (2.11)$$

$$C_1 := \frac{p(1)}{g_1^2 + v_1^2} \sqrt{\sum_{n=1}^{\infty} \frac{1}{\lambda_n^2} \left| v_1 \phi_n(1) - g_1 \frac{d\phi_n}{dz}(1) \right|^2} = \frac{1}{\sqrt{g_1^2 + v_1^2}} \|\bar{x}\|_{2,r}, \quad (2.12)$$

$$\tilde{C} := \sqrt{\sum_{n=1}^{\infty} \frac{1}{\lambda_n^2} \left( \int_0^1 r(z) |\phi_n(z)| dz \right)^2}, \quad (2.13)$$

$$\|x[t]\|_{2,r} := \left( \int_0^1 r(z) |x(t,z)|^2 dz \right)^{1/2}, \quad (2.14)$$

$\tilde{x} \in C^2([0,1])$ *is the unique solution of the boundary value problem* $\frac{d}{dz}\left( p(z) \frac{d\tilde{x}}{dz}(z) \right) - q(z)\tilde{x}(z) = 0$ *for* $z \in [0,1]$ *with* $g_0 \tilde{x}(0) + v_0 \frac{d\tilde{x}}{dz}(0) = \sqrt{g_0^2 + v_0^2}$, $g_1 \tilde{x}(1) + v_1 \frac{d\tilde{x}}{dz}(1) = 0$ *and* $\bar{x} \in C^2([0,1])$ *is the unique solution of the boundary value problem* $\frac{d}{dz}\left( p(z) \frac{d\bar{x}}{dz}(z) \right) - q(z)\bar{x}(z) = 0$ *for* $z \in [0,1]$ *with* $g_0 \bar{x}(0) + v_0 \frac{d\bar{x}}{dz}(0) = 0$ *and* $g_1 \bar{x}(1) + v_1 \frac{d\bar{x}}{dz}(1) = \sqrt{g_1^2 + v_1^2}$.

Finally, the third main result of the section provides an ISS estimate in a weighted $L^1$ norm for a special case of the parabolic PDE (2.6) with Dirichlet boundary conditions.

**Theorem 2.4:** *Let $a > 0$, $b, k \in \Re$ be given constants with $k < a\pi^2 + \frac{b^2}{4a}$. For every $x_0 \in C^2([0,1])$, $v \in C^1(\Re_+ \times [0,1])$, let $\Phi(x_0, v) \subseteq C^1(\Re_+; \Re^2)$ be the non-empty set for which the evolution equation*

$$\frac{\partial x}{\partial t}(t,z) = a \frac{\partial^2 x}{\partial z^2}(t,z) + b \frac{\partial x}{\partial z}(t,z) + kx(t,z) + v(t,z) \text{ for } (t,z) \in (0,+\infty) \times (0,1) \quad (2.15)$$

*with*

$$x(t,0) = d_0(t), \quad x(t,1) = d_1(t) \text{ for } t \geq 0 \quad (2.16)$$

$(d_0, d_1) \in \Phi(x_0, v)$ *and initial condition* $x_0 \in C^2([0,1])$ *has a unique solution* $x \in C^0(\Re_+ \times [0,1]) \cap C^1((0,+\infty) \times [0,1])$ *with* $x[t] \in C^2([0,1])$ *for all* $t \geq 0$ *and* $x(0,z) = x_0(z)$ *for all* $z \in [0,1]$. *Then for every* $x_0 \in C^2([0,1])$, $v \in C^1(\Re_+ \times [0,1])$ *and* $(d_0, d_1) \in \Phi(x_0, v)$, *the unique solution* $x \in C^0(\Re_+ \times [0,1]) \cap C^1((0,+\infty) \times [0,1])$ *of the evolution equation (2.15) with (2.16) and initial condition* $x_0 \in C^2([0,1])$ *satisfies the following estimate for all* $t \geq 0$:

$$\|x[t]\|_{1,w} \leq \exp\left(-\left(a\pi^2 - k + \frac{b^2}{4a}\right)t\right)\|x[0]\|_{1,w} + \frac{4a^2\pi}{4a^2\pi^2 + b^2 - 4ak}\left(\max_{0\leq s\leq t}(|d_0(s)|) + \exp\left(\frac{b}{2a}\right)\max_{0\leq s\leq t}(|d_0(s)|)\right)$$

$$+ \frac{4a}{4a^2\pi^2 + b^2 - 4ak}\max_{0\leq s\leq t}\left(\max_{0\leq z\leq 1}\left(|v(s,z)|\exp\left(\frac{bz}{2a}\right)\sin(\pi z)\right)\right)$$ (2.17)

*where*

$$\|x[t]\|_{1,w} := \int_0^1 \exp\left(\frac{bz}{2a}\right)\sin(\pi z)|x(t,z)|dz .$$ (2.18)

**Remarks on the three main results:**
  (a) The ISS estimate (2.8) is different from the ISS estimates (2.10) and (2.17): estimate (2.8) is expressed by means of max-operators (the max-formulation of the ISS property; see [12]), while the ISS estimates (2.10), (2.17) are expressed by means of summation.
  (b) The proofs of Theorem 2.2 and Theorem 2.4 are similar: they are both exploiting a finite-difference scheme for the approximation of the solution of the PDE. On the other hand, the proof of Theorem 2.3 relies on an eigenfunction expansion, which was also utilized in [14].
  (c) Theorem 2.4 is only applicable to the PDE (2.15) with Dirichlet boundary conditions (2.16), while Theorem 2.2 and Theorem 2.3 are applicable to much more general parabolic PDEs with various boundary conditions (Dirichlet, Neumann or Robin). Assumptions (H1), (H2) are not restrictive. To see this, notice that the boundary condition $x(t,0) = d_0(t)$ (which violates Assumption (H1)) can always be written as $-x(t,0) = \tilde{d}_0(t)$, where $\tilde{d}_0(t) \equiv -d_0(t)$ and thus Assumption (H1) holds with $g_0 = -1$. Similar manipulations may be done for the satisfaction of Assumption (H2) as well (if needed).
  (d) Assumption (H4) requires the existence of a positive eigenfunction of the Sturm-Liouville operator $A: D \to L_r^2([0,1])$ defined by (2.1). The problem of the existence of positive eigenfunctions for elliptic operators has been studied in the literature (see [25] on page 112 and references therein). However, it must be noted that there is a degree of freedom in the selection of the boundary conditions: no specific boundary conditions are assumed to hold for the function $\eta \in C^2([0,1];(0,+\infty))$: the function is only required to satisfy the (boundary) inequalities $g_0\eta(0) + v_0\eta'(0) < 0$ and $g_1\eta(1) + v_1\eta'(1) > 0$.
  (e) Since the equilibrium points that correspond to the constant disturbances $d_0(t) \equiv \sqrt{g_0^2 + v_0^2}$ and $d_1(t) \equiv 0$, $u(t,z) \equiv 0$, or $d_1(t) \equiv \sqrt{g_1^2 + v_1^2}$ and $d_0(t) \equiv 0$, $u(t,z) \equiv 0$, are the functions $\tilde{x} \in C^2([0,1])$ and $\bar{x} \in C^2([0,1])$, respectively, it follows that the gains $\gamma_0, \gamma_1 > 0$ of the boundary inputs that are involved to the ISS estimate
  $$\|x[t]\|_{2,r} \leq M\exp(-\sigma t)\|x[0]\|_{2,r} + \gamma_0 \max_{0\leq s\leq t}(|d_0(s)|) + \gamma_1 \max_{0\leq s\leq t}(|d_1(s)|), \text{ for all } t\geq 0$$
  for certain constants $M, \sigma$, must satisfy the inequalities
  $$\gamma_0 \geq C_0, \ \gamma_1 \geq C_1$$
  where $C_0, C_1 > 0$ are given by (2.11) and (2.12). On the other hand, Theorem 2.3 guarantees that
  $$\gamma_0 \leq \sqrt{(1+\varepsilon^{-1})(1+\omega)}\, C_0, \ \gamma_1 \leq \sqrt{(1+\varepsilon)(1+\omega)}\, C_1, \text{ for all } \varepsilon, \omega > 0$$
  Consequently, we can guarantee that the estimation of the gain made by Theorem 2.3 is sharp. Moreover, formulas (2.11), (2.12) guarantee that the gains of the boundary disturbances can be computed **without exact knowledge** of the eigenvalues and the eigenfunctions of the Sturm-Liouville operator $A: D \to L_r^2([0,1])$ defined by (2.1), (2.2). The only thing we need to know about the eigenvalues and the eigenfunctions of the Sturm-Liouville operator $A: D \to L_r^2([0,1])$ defined by (2.1), (2.2) is that Assumptions (H1), (H2), (H3) hold.

In order to understand how well the gains of the inputs $d_0, d_1 \in$ are estimated by the inequalities (2.8), (2.10) and (2.17), we next specialize the results to the heat equation with Dirichlet boundary conditions.

**Corollary 2.5:** *Let $a > 0$ be a given constant. For every $x_0 \in C^2([0,1])$ let $\Phi(x_0) \subseteq C^0(\Re_+; \Re^2)$ be the non-empty set for which the evolution equation*

$$\frac{\partial x}{\partial t}(t,z) = a \frac{\partial^2 x}{\partial z^2}(t,z) \quad \text{for } (t,z) \in (0,+\infty) \times (0,1) \tag{2.19}$$

*with (2.16), $(d_0, d_1) \in \Phi(x_0)$ and initial condition $x_0 \in C^2([0,1])$ has a unique solution $x \in C^0(\Re_+ \times [0,1]) \cap C^1((0,+\infty) \times [0,1])$ with $x[t] \in C^2([0,1])$ for all $t \geq 0$ and $x(0,z) = x_0(z)$ for all $z \in [0,1]$. Then for every $x_0 \in C^2([0,1])$ and $(d_0, d_1) \in \Phi(x_0)$, the unique solution $x \in C^0(\Re_+ \times [0,1]) \cap C^1((0,+\infty) \times [0,1])$ of the evolution equation (2.19) with (2.16) and initial condition $x_0 \in C^2([0,1])$ satisfies the following estimates for all $t \geq 0$:*

$$\int_0^1 \sin(\pi z)|x(t,z)|dz \leq \exp(-a\pi^2 t)\int_0^1 \sin(\pi z)|x(0,z)|dz + \frac{1}{\pi}\left(\max_{0 \leq s \leq t}(|d_0(s)|) + \max_{0 \leq s \leq t}(|d_0(s)|)\right) \tag{2.20}$$

$$\sqrt{\int_0^1 x^2(t,z)dz} \leq \sqrt{\frac{\exp(-a\pi^2 t)}{2 - \exp(-a\pi^2 t)}}\sqrt{\int_0^1 x^2(0,z)dz} + \frac{1}{\sqrt{3}}\left(\max_{0 \leq s \leq t}(|d_0(s)|) + \max_{0 \leq s \leq t}(|d_1(s)|)\right), \tag{2.21}$$

$$\max_{0 \leq z \leq 1}\left(\frac{\sin(\theta+\varphi)|x(t,z)|}{\sin(\theta+z\varphi)}\right) \leq \max\left(\exp(-\sigma t)\max_{0 \leq z \leq 1}\left(\frac{\sin(\theta+\varphi)|x(0,z)|}{\sin(\theta+z\varphi)}\right), \frac{\sin(\theta+\varphi)}{\sin(\theta)}\max_{0 \leq s \leq t}(|d_0(s)|), \max_{0 \leq s \leq t}(|d_1(s)|)\right),$$

$$\text{for all } \sigma \in (0, a\pi^2), \; \theta \in (0, \pi - \varphi) \text{ with } \varphi := \sqrt{\frac{\sigma}{a}} \tag{2.22}$$

**Remark 2.6:**
(a) It is clear that the estimation from (2.20), (2.21) and (2.22) of the gains of the inputs $d_0, d_1$ is sharp, since the constant disturbances $d_0(t) \equiv 1$ or $d_1(t) \equiv 1$ lead to equilibrium profiles for which the corresponding weighted norms are exactly equal to the estimated gains.

(b) By selecting $\sigma = a(\pi - 2\theta)^2$, we get from (2.22) the following estimate, which holds for all $\theta \in (0, \pi/2)$ and $t \geq 0$:

$$\max_{0 \leq z \leq 1}\left(\frac{|x(t,z)|}{\sin(\theta + z(\pi - 2\theta))}\right) \leq \max\left(\exp(-a(\pi - 2\theta)^2 t)\max_{0 \leq z \leq 1}\left(\frac{|x(0,z)|}{\sin(\theta + z(\pi - 2\theta))}\right), \frac{\max_{0 \leq s \leq t}(|d_0(s)|)}{\sin(\theta)}, \frac{\max_{0 \leq s \leq t}(|d_1(s)|)}{\sin(\theta)}\right) \tag{2.23}$$

Estimate (2.23) is a refinement of the well-known maximum principle for the heat equation (see [11] on page 215): the classical maximum principle coincides with estimate (2.23) for $\theta = \pi/2$. Taking $\theta = \pi/4$, we obtain from (2.23) the following simple estimate

$$\max_{0 \leq z \leq 1}(|x(t,z)|) \leq \sqrt{2}\max\left(\exp\left(-\frac{a\pi^2 t}{4}\right)\max_{0 \leq z \leq 1}(|x(0,z)|), \max_{0 \leq s \leq t}(|d_0(s)|), \max_{0 \leq s \leq t}(|d_1(s)|)\right) \tag{2.24}$$

which holds for all $t \geq 0$.

# 3. Application to Thermoelasticity Models Via Small-Gain Argument

Quasi-static models of thermoelasticity in one spatial dimension deal with the PDE (2.19) with boundary conditions

$$x(t,0) = \int_0^1 g_0(z)x(t,z)dz \text{ and } x(t,1) = \int_0^1 g_1(z)x(t,z)dz \qquad (3.1)$$

where $g_0, g_1 \in C^1([0,1])$ are given functions (see [7,8]). The variable $x(t,z)$ is the dimensionless entropy at position $z \in [0,1]$ (dimensionless spatial variable) and time $t \geq 0$ (dimensionless time variable) of a slab of homogeneous and isotropic material. Conditions for the functions $g_0, g_1 \in C^1([0,1])$ can be imposed (see [23]) so that there exists a non-empty set

$$\Delta \subseteq \left\{ x \in C^2([0,1]) : x(0) - \int_0^1 g_0(z)x(z)dz = x(1) - \int_0^1 g_1(z)x(z)dz = 0 \right\}$$

which is dense in $L^2$ and such that for every $x_0 \in \Delta$ the evolution equation (2.19) with (3.1) has a unique solution $x \in C^0(\Re_+ \times [0,1]) \cap C^1((0,+\infty) \times [0,1])$ with $x[t] \in C^2([0,1])$ for all $t \geq 0$ and $x(0,z) = x_0(z)$ for all $z \in [0,1]$. However, we will not deal here with the existence/uniqueness problem for the solution of (2.19) with (3.1).

The present section is devoted to the stability analysis of the zero solution for the dynamical system described by (2.19) with (3.1). The following results use the ISS estimates of the previous sections in conjunction with small-gain arguments.

**Theorem 3.1:** *Suppose that there exist constants $\varphi \in (0, \pi)$, $\theta \in (0, \pi - \varphi)$ such that*

$$\int_0^1 |g_0(z)| \sin(\theta + z\varphi)dz < \sin(\theta) \text{ and } \int_0^1 |g_1(z)| \sin(\theta + z\varphi)dz < \sin(\theta + \varphi) \qquad (3.2)$$

*Then there exist constants $M, \delta > 0$ such that for every $x_0 \in \Delta$ the unique solution $x \in C^0(\Re_+ \times [0,1]) \cap C^1((0,+\infty) \times [0,1])$ of the evolution equation (2.19) with (3.1) with $x[t] \in C^2([0,1])$ for all $t \geq 0$ and $x(0,z) = x_0(z)$ for all $z \in [0,1]$ satisfies the following estimate for all $t \geq 0$:*

$$\|x[t]\|_\infty \leq M \exp(-\delta t) \|x[0]\|_\infty \qquad (3.3)$$

*where $\|x[t]\|_\infty := \max_{0 \leq z \leq 1}(|x(t,z)|)$.*

**Theorem 3.2:** *Suppose that*

$$\left( \int_0^1 g_0^2(z)dz \right)^{1/2} + \left( \int_0^1 g_1^2(z)dz \right)^{1/2} < \sqrt{3} \qquad (3.4)$$

*Then there exist constants $M, \delta > 0$ such that for every $x_0 \in \Delta$ the unique solution $x \in C^0(\Re_+ \times [0,1]) \cap C^1((0,+\infty) \times [0,1])$ of the evolution equation (2.19) with (3.1) with $x[t] \in C^2([0,1])$ for all $t \geq 0$ and $x(0,z) = x_0(z)$ for all $z \in [0,1]$ satisfies the following estimate for all $t \geq 0$:*

$$\|x[t]\|_2 \leq M \exp(-\delta t)\|x[0]\|_2 \qquad (3.5)$$

where $\|x[t]\|_2 := \left(\int_0^1 |x(t,z)|^2 dz\right)^{1/2}$.

**Theorem 3.3:** *Suppose that* $g_0(0) = g_0(1) = g_1(0) = g_1(1) = 0$ *and that*

$$\sup_{0<z<1}\left(\frac{|g_0(z)|}{\sin(\pi z)}\right) + \sup_{0<z<1}\left(\frac{|g_1(z)|}{\sin(\pi z)}\right) < \pi \qquad (3.6)$$

*Then there exist constants* $M, \delta > 0$ *such that for every* $x_0 \in \Delta$ *the unique solution* $x \in C^0(\Re_+ \times [0,1]) \cap C^1((0,+\infty) \times [0,1])$ *of the evolution equation (2.19) with (3.1) with* $x[t] \in C^2([0,1])$ *for all* $t \geq 0$ *and* $x(0,z) = x_0(z)$ *for all* $z \in [0,1]$ *satisfies the following estimate for all* $t \geq 0$:

$$\|x[t]\|_{1,w} \leq M \exp(-\delta t)\|x[0]\|_{1,w} \qquad (3.7)$$

where $\|x[t]\|_{1,w} := \int_0^1 \sin(\pi z)|x(t,z)|dz$.

**Remark 3.4:**
  (a) The proofs of Theorem 3.1, Theorem 3.2 and Theorem 3.3 show that conditions (3.4), (3.6) and the pair of conditions (3.2) are small-gain conditions for the heat equation (2.19) with a pair of boundary feeback inputs given in (3.1). Similar stability conditions have been used in [7,8] and in [9,10,15] (for numerical stability of discretization schemes). It should be noticed that conditions (3.2) provide a generalization of the stability conditions given in [7,8]

$$\int_0^1 |g_0(z)|dz < 1 \text{ and } \int_0^1 |g_1(z)|dz < 1 \qquad (3.8)$$

  in the sense that if conditions (3.8) hold then there exist constants $\varphi \in (0,\pi)$, $\theta \in (0,\pi-\varphi)$ such that inequalities (3.2) hold. This can be guaranteed by continuity of the functions $P_0(\theta,\varphi) := \int_0^1 |g_0(z)|\sin(\theta+z\varphi)dz - \sin(\theta)$, $P_1(\theta,\varphi) := \int_0^1 |g_1(z)|\sin(\theta+z\varphi)dz - \sin(\theta+\varphi)$ and the fact that inequalities (3.8) imply the inequalities $P_0\left(\frac{\pi}{2},0\right) < 0$, $P_1\left(\frac{\pi}{2},0\right) < 0$.

  (b) There is a difference between our results and the results in [7,8]. The results in [7,8] guarantee that the mapping $t \to \|x[t]\|_\infty$ is decreasing under conditions (3.6), while our results guarantee the exponential convergence of $\|x[t]\|_\infty$ or $\|x[t]\|_2$ or $\|x[t]\|_{1,w}$ to zero under conditions (3.2), (3.4) or (3.6), respectively. However, our results do not exclude the possibility that $\|x[t]\|_\infty$ or $\|x[t]\|_2$ or $\|x[t]\|_{1,w}$ might increase for a transient period.

## 4. Proofs of Main Results

We start with the proof of Theorem 2.1.

**Proof of Theorem 2.1:** Let arbitrary $d_0, d_1 \in C^2(\Re_+)$, $u \in C^1(\Re_+ \times [0,1])$ and $x_0 \in X(d_0(0), d_1(0))$ be given. For the existence/uniqueness of the evolution equation (2.6) with (2.7) and initial condition $x_0 \in X(d_0(0), d_1(0))$, we simply apply the transformation $x(t, z) = y(t, z) + p_0(z) d_0(t) + d_1(t) p_1(z)$, where $y \in C^0(\Re_+ \times [0,1]; \Re) \cap C^1((0, +\infty) \times [0,1]; \Re)$ is the unique function that satisfies $y[t] \in X(0,0)$ (recall (2.5)) for all $t \geq 0$, $y(0, z) = x(0, z) - p_0(z) d_0(0) - p_1(z) d_1(0)$ for all $z \in [0,1]$ and

$$\frac{\partial y}{\partial t}(t, z) + (Ay[t])(z) = -\dot{d}_0(t) p_0(z) - \dot{d}_1(t) p_1(z) - d_0(t)(Ap_0)(z) - d_1(t)(Ap_1)(z) + u(t, z),$$

for all $(t, z) \in (0, +\infty) \times (0,1)$ (4.1)

where $p_0(z) := (c_0 + c_1 z)(1-z)^2$, $p_1(z) := (C_0 + C_1(1-z))z^2$, $c_0 = \dfrac{g_0 - 2v_0}{(g_0 - 2v_0)^2 + v_0^2}$, $c_1 = \dfrac{v_0}{(g_0 - 2v_0)^2 + v_0^2}$, $C_0 = \dfrac{g_1 + 2v_1}{(g_1 + 2v_1)^2 + v_1^2}$, $C_1 = \dfrac{-v_1}{(g_1 + 2v_1)^2 + v_1^2}$ (notice that Assumptions (H1), (H2) guarantee that $(g_0 - 2v_0)^2 + v_0^2 > 0$, $(g_1 + 2v_1)^2 + v_1^2 > 0$). The existence/uniqueness of $y \in C^0(\Re_+ \times [0,1]; \Re) \cap C^1((0, +\infty) \times [0,1]; \Re)$ that satisfies $y[t] \in X(0,0)$ for all $t \geq 0$, $y(0, z) = x(0, z) - p_0(z) d_0(0) - p_1(z) d_1(0)$ for all $z \in [0,1]$ and (4.1) is guaranteed by Theorem 3.1 in [14]. The proof is complete. ◁

The proof of Theorem 2.2 utilizes the following technical lemma. Its proof is provided in the Appendix.

**Lemma 4.1:** *Let* $\{y(j) \geq 0\}_{j=0}^\infty$, $\{\varphi(j) \geq 0\}_{j=0}^\infty$ *be sequences that satisfy the following properties:*

**(P1)** $j \in \{0,1,\ldots,m\}, y(j) \leq K\varphi(j) \Rightarrow \varphi(j+1) \leq a\varphi(j) + g$

**(P2)** $j \in \{0,1,\ldots,m\}, y(j) \geq K\varphi(j) \Rightarrow \varphi(j+1) \leq \beta y(j) + g$

*for certain real constants* $K, g, \beta > 0$, $a \in (0,1)$ *and for certain positive integer* $m > 0$. *Then the following inequality holds for all* $j \in \{0,1,\ldots,m\}$:

$$\varphi(j) \leq \max\left( \max(K^{-1}, \beta) \max_{s=0,\ldots,j}(y(s)), a^j \varphi(0) \right) + \frac{g}{1-a} \quad (4.2)$$

**Proof of Theorem 2.2:** Let $\eta \in C^2([0,1]; (0, +\infty))$ be the function involved in Assumption (H4). Notice that we may assume that $\eta \in C^2([-1, 2])$ (we can use an arbitrary extension of the function $\eta$ out of the interval $[0,1]$).

In the following exposition, we will use the notation

$$a(z) := \frac{p(z)}{r(z)}, \quad b(z) := \frac{p'(z)}{r(z)}, \quad c(z) := -\frac{q(z)}{r(z)}, \text{ for } z \in [0,1] \quad (4.3)$$

It suffices to show that there exists a constant $\Gamma > 0$ such that for every time instants $0 < t_0 \leq T$ and for every $\varepsilon > 0$, the following estimate holds for sufficiently large $N > 1$:

$$V_N(T) \leq \max\left( \max\left( \frac{\|d_0\|}{\gamma_0(N^{-1})}, \frac{\|d_1\|}{\gamma_1(N^{-1})} \right), \exp\left( -\sigma\left( 1 - \max_{0 \leq z \leq 1}(\omega(z, N^{-1})) \right)(T - t_0) \right) V_N(t_0) \right) + \frac{\Gamma \varepsilon + \|u\|}{\sigma\left( 1 - \max_{0 \leq z \leq 1}(\omega(z, N^{-1})) \right)},$$
(4.4)

where $\sigma > 0$ is the constant involved in Assumption (H4), $\|d_0\| := \max_{0 \leq s \leq T}(|d_0(s)|)$, $\|d_1\| := \max_{0 \leq s \leq T}(|d_1(s)|)$,

$\|u\| := \max_{0 \leq s \leq T}\left( \max_{0 \leq z \leq 1}\left( \frac{|u(s,z)|}{\eta(z)} \right) \right)$,

$$V_N(t) := \max_{i=1,\ldots N-1}\left( \frac{|x(t, iN^{-1})|}{\eta(iN^{-1})} \right), \tag{4.5}$$

$$\gamma_0(h) := -g_0 \eta(0) - v_0 \eta'(0) + v_0 \left( \eta'(0) - \frac{4\eta(h) - 3\eta(0) - \eta(2h)}{2h} \right), \text{ for } h \in (0,1) \tag{4.6}$$

$$\gamma_1(h) := g_1 \eta(1) + v_1 \eta'(1) + v_1 \left( \frac{3\eta(1) - 4\eta(1-h) + \eta(1-2h)}{2h} - \eta'(1) \right), \text{ for } h \in (0,1) \tag{4.7}$$

and the function $\omega : [0,1] \times [0,1] \to \Re$ is given by the formulas

$$\omega(z, 0) := 0, \text{ for all } z \in [0,1] \tag{4.8}$$

$$\omega(z, h) := \sigma^{-1}\left( a(z) + \frac{h}{2}b(z) \right) \frac{\eta(z+h) - \eta(z) - h\eta'(z) - \frac{h^2}{2}\eta''(z)}{h^2 \eta(z)}$$
$$+ \sigma^{-1}\left( a(z) - \frac{h}{2}b(z) \right) \frac{\eta(z-h) - \eta(z) + h\eta'(z) - \frac{h^2}{2}\eta''(z)}{h^2 \eta(z)}, \text{ for all } (z,h) \in [0,1] \times (0,1] \tag{4.9}$$

Indeed, since

$$\max_{0 \leq z \leq 1}\left( \frac{|x(T,z)|}{\eta(z)} \right) \leq V_N(t) + N^{-1} \max_{0 \leq z \leq 1}\left( \left| \frac{\partial}{\partial z}\left( \frac{x(T,z)}{\eta(z)} \right) \right| \right) \tag{4.10}$$

we obtain from (4.4) that

$$\max_{0 \leq z \leq 1}\left( \frac{|x(T,z)|}{\eta(z)} \right) \leq \max\left( \max\left( \frac{\|d_0\|}{\gamma_0(N^{-1})}, \frac{\|d_1\|}{\gamma_1(N^{-1})} \right), \exp\left( -\sigma\left( 1 - \max_{0 \leq z \leq 1}(\omega(z, N^{-1})) \right)(T - t_0) \right) \max_{0 \leq z \leq 1}\left( \frac{|x(t_0, z)|}{\eta(z)} \right) \right)$$
$$+ N^{-1} M + \frac{\Gamma \varepsilon + \|u\|}{\sigma\left( 1 - \max_{0 \leq z \leq 1}(\omega(z, N^{-1})) \right)}$$

where $M := \max_{0 \leq z \leq 1}\left( \left| \frac{\partial}{\partial z}\left( \frac{x(T,z)}{\eta(z)} \right) \right| \right)$. Since the function $\omega : [0,1] \times [0,1] \to \Re$ given by (4.8), (4.9) is a continuous function (and thus uniformly continuous on the compact set $[0,1] \times [0,1]$) and since the limits $\lim_{h \to 0^+}(\gamma_0(h))$, $\lim_{h \to 0^+}(\gamma_1(h))$ exist (see definitions (4.6), (4.7)), by letting $N \to +\infty$, we get (using (4.8)) for all $\varepsilon > 0$:

$$\max_{0 \leq z \leq 1}\left( \frac{|x(T,z)|}{\eta(z)} \right) \leq \max\left( \max\left( \frac{\|d_0\|}{-g_0\eta(0) - v_0\eta'(0)}, \frac{\|d_1\|}{g_1\eta(1) - v_1\eta'(1)} \right), \exp(-\sigma(T - t_0)) \max_{0 \leq z \leq 1}\left( \frac{|x(t_0, z)|}{\eta(z)} \right) \right) + \frac{\Gamma \varepsilon + \|u\|}{\sigma}$$

Since the above inequality holds for all $t_0 \in (0,T]$, $\varepsilon > 0$ and since $\lim\limits_{t_0 \to 0^+}\left(\max\limits_{0 \leq z \leq 1}\left(\frac{|x(t_0,z)|}{\eta(z)}\right)\right) = \max\limits_{0 \leq z \leq 1}\left(\frac{|x_0(z)|}{\eta(z)}\right)$, the above inequality implies that (2.8) holds for all $t > 0$. Inequality (2.8) holds for $t = 0$ as well.

Therefore, we have to prove that there exists a constant $\Gamma > 0$ such that for every time instants $0 < t_0 \leq T$ and for every $\varepsilon > 0$, estimate (4.4) holds for sufficiently large $N > 1$.

Let $0 < t_0 \leq T$ be arbitrary time instants. For every $h \in (0,1)$, $\delta > 0$, $t \geq t_0$, $z \in [h, 1-h]$ we obtain:

$$x(t+\delta, z) = x(t,z) + \delta \frac{\partial x}{\partial t}(t,z) + \int_t^{t+\delta}\left(\frac{\partial x}{\partial t}(s,z) - \frac{\partial x}{\partial t}(t,z)\right)ds \quad (4.11)$$

$$x(t, z+h) - 2x(t,z) + x(t,z-h) = h^2 \frac{\partial^2 x}{\partial z^2}(t,z)$$
$$+ \int_z^{z+h}\int_z^s\left(\frac{\partial^2 x}{\partial z^2}(t,w) - \frac{\partial^2 x}{\partial z^2}(t,z)\right)dw\,ds + \int_z^{z-h}\int_z^s\left(\frac{\partial^2 x}{\partial z^2}(t,w) - \frac{\partial^2 x}{\partial z^2}(t,z)\right)dw\,ds \quad (4.12)$$

$$x(t,z+h) - x(t,z-h) = 2h\frac{\partial x}{\partial z}(t,z)$$
$$+ \int_z^{z+h}\int_z^s\left(\frac{\partial^2 x}{\partial z^2}(t,w) - \frac{\partial^2 x}{\partial z^2}(t,z)\right)dw\,ds - \int_z^{z-h}\int_z^s\left(\frac{\partial^2 x}{\partial z^2}(t,w) - \frac{\partial^2 x}{\partial z^2}(t,z)\right)dw\,ds \quad (4.13)$$

Let $h^* \in (0,1)$ be a constant so that the following inequalities hold for all $h \in (0, h^*)$:

$$3v_0 - 2hg_0 > 0, \quad 3v_1 + 2hg_1 > 0, \quad \frac{v_0}{3v_0 - 2hg_0} \leq 1, \quad \frac{v_1}{3v_1 + 2hg_1} \leq 1 \quad (4.14)$$

$$\min_{0 \leq z \leq 1}(a(z)) > \frac{h}{2}\max_{0 \leq z \leq 1}(|b(z)|) \quad (4.15)$$

$$a(h) + \frac{h}{2}b(h) > \frac{v_0}{3v_0 - 2hg_0}\left(a(h) - \frac{h}{2}b(h)\right), \quad a(1-h) - \frac{h}{2}b(1-h) > \frac{v_1}{3v_1 + 2hg_1}\left(a(1-h) + \frac{h}{2}b(1-h)\right) \quad (4.16)$$

$$-v_0\left(\frac{4\eta(h) - 3\eta(0) - \eta(2h)}{2h} - \eta'(0)\right) > g_0\eta(0) + v_0\eta'(0), \quad v_1\left(\frac{3\eta(1) - 4\eta(1-h) + \eta(1-2h)}{2h} - \eta'(1)\right) + g_1\eta(1) + v_1\eta'(1) > 0 \quad (4.17)$$

Notice that the existence of $h^* \in (0,1)$ so that (4.14) holds for all $h \in (0, h^*)$ is a direct consequence of Assumptions (H1), (H2). The existence of $h^* \in (0,1)$ so that (4.17) holds for all $h \in (0, h^*)$ is a direct consequence of Assumption (H4) and the facts that $g_0\eta(0) + v_0\eta'(0) < 0$ and $g_1\eta(1) + v_1\eta'(1) > 0$. Finally, (4.15) and (4.16) is a direct consequence of the continuity of the functions $a(z) = \frac{p(z)}{r(z)}$, $b(z) = \frac{p'(z)}{r(z)}$ and the fact that $a(z) > 0$ for all $z \in [0,1]$.

We next select an integer $N > \frac{1}{h^*}$ and define:

$h := N^{-1}$, $\delta := \lambda h^2$, $x_i(t) := x(t, ih)$, $u_i(t) := u(t, ih)$, $a_i := a(ih)$, $b_i := b(ih)$, $c_i := c(ih)$, $\eta_i := \eta(ih)$ $i = 0,...,N$, (4.18)

where

$$\lambda \in \left(0, \frac{1}{1 + 2\bar{a} + C}\right) \quad (4.19)$$

is a constant to be selected and

$$\bar{a} := \max_{0 \le z \le 1}(a(z)), \quad C := \min_{0 \le z \le 1}(|c(z)|) \tag{4.20}$$

Notice that (4.19) in conjunction with definitions (4.18), (4.20) and (4.14), (4.15) guarantee the following inequalities for all $N > \frac{1}{h^*}$:

$$1 - 2\lambda a_i + \lambda h^2 c_i > 0, \text{ for } i = 0,\ldots,N \tag{4.21}$$

$$1 - 2\lambda a_1 + \lambda h^2 c_1 + \lambda \frac{4v_0}{3v_0 - 2hg_0}\left(a_1 - \frac{h}{2}b_1\right) > 0 \tag{4.22}$$

$$1 - 2\lambda a_{N-1} + \lambda h^2 c_{N-1} + \lambda\left(a_{N-1} + \frac{h}{2}b_{N-1}\right)\frac{4v_1}{3v_1 + 2hg_1} > 0 \tag{4.23}$$

It follows from (2.6), (4.3), (4.11), (4.12), (4.13) and the fact $\delta := \lambda h^2$ that for every $h \in (0,1)$, $\lambda \in \left(0, \frac{1}{1+2\bar{a}+C}\right)$, $t \ge t_0$, $z \in [h, 1-h]$, the following equality holds:

$$x(t+\delta, z) = \left(1 - 2\lambda a(z) + \lambda h^2 c(z)\right)x(t,z) + \lambda\left(a(z) + \frac{h}{2}b(z)\right)x(t, z+h)$$
$$+ \lambda\left(a(z) - \frac{h}{2}b(z)\right)x(t, z-h) + \lambda h^2 u(t,z) + P(\delta, h, t, z) \tag{4.24}$$

where

$$P(\delta,h,t,z) := \int_t^{t+\delta}\left(\frac{\partial x}{\partial t}(s,z) - \frac{\partial x}{\partial t}(t,z)\right)ds - \lambda\left(a(z) + \frac{h}{2}b(z)\right)\int_z^{z+h}\int_z^s\left(\frac{\partial^2 x}{\partial z^2}(t,w) - \frac{\partial^2 x}{\partial z^2}(t,z)\right)dw\,ds$$
$$- \lambda\left(a(z) - \frac{h}{2}b(z)\right)\int_z^{z-h}\int_z^s\left(\frac{\partial^2 x}{\partial z^2}(t,w) - \frac{\partial^2 x}{\partial z^2}(t,z)\right)dw\,ds \tag{4.25}$$

Furthermore, it follows from (4.24) and definitions (4.18) that for every $t \ge t_0$, the following equalities hold:

$$x_i(t+\delta) = \left(1 - 2\lambda a_i + \lambda h^2 c_i\right)x_i(t) + \lambda\left(a_i + \frac{h}{2}b_i\right)x_{i+1}(t) + \lambda\left(a_i - \frac{h}{2}b_i\right)x_{i-1}(t) + \lambda h^2 u_i(t) + P(\delta, h, t, ih),$$
$$\text{for } i = 1,\ldots,N-1 \tag{4.26}$$

Using (2.7) and the following equalities

$$x_1(t) - x_0(t) - \int_0^h\int_0^s\left(\frac{\partial^2 x}{\partial z^2}(t,w) - \frac{\partial^2 x}{\partial z^2}(t,0)\right)dw\,ds = h\frac{\partial x}{\partial z}(t,0) + \frac{h^2}{2}\frac{\partial^2 x}{\partial z^2}(t,0)$$

$$x_2(t) - x_0(t) - \int_0^{2h}\int_0^s\left(\frac{\partial^2 x}{\partial z^2}(t,w) - \frac{\partial^2 x}{\partial z^2}(t,0)\right)dw\,ds = 2h\frac{\partial x}{\partial z}(t,0) + 2h^2\frac{\partial^2 x}{\partial z^2}(t,0)$$

$$x_{N-1}(t) - x_N(t) - \int_{1-h}^1\int_s^1\left(\frac{\partial^2 x}{\partial z^2}(t,w) - \frac{\partial^2 x}{\partial z^2}(t,1)\right)dw\,ds = -h\frac{\partial x}{\partial z}(t,1) + \frac{h^2}{2}\frac{\partial^2 x}{\partial z^2}(t,1)$$

$$x_{N-2}(t) - x_N(t) - \int_{1-2h}^1\int_s^1\left(\frac{\partial^2 x}{\partial z^2}(t,w) - \frac{\partial^2 x}{\partial z^2}(t,1)\right)dw\,ds = -2h\frac{\partial x}{\partial z}(t,1) + 2h^2\frac{\partial^2 x}{\partial z^2}(t,1) \tag{4.27}$$

we obtain the equalities

$$(3v_0 - 2hg_0)x_0(t) = -2hd_0(t) + 4v_0 x_1(t) - v_0 x_2(t)$$
$$+ v_0 \left( \int_0^{2h}\!\!\int_0^s \left( \frac{\partial^2 x}{\partial z^2}(t,w) - \frac{\partial^2 x}{\partial z^2}(t,0) \right) dw\,ds - 4 \int_0^h\!\!\int_0^s \left( \frac{\partial^2 x}{\partial z^2}(t,w) - \frac{\partial^2 x}{\partial z^2}(t,0) \right) dw\,ds \right) \quad (4.28)$$

and

$$(3v_1 + 2hg_1)x_N(t) = 2hd_1(t) + 4v_1 x_{N-1}(t) - v_1 x_{N-2}(t)$$
$$+ v_1 \left( \int_{1-2h}^{1}\!\!\int_s^1 \left( \frac{\partial^2 x}{\partial z^2}(t,w) - \frac{\partial^2 x}{\partial z^2}(t,1) \right) dw\,ds - 4 \int_{1-h}^{1}\!\!\int_s^1 \left( \frac{\partial^2 x}{\partial z^2}(t,w) - \frac{\partial^2 x}{\partial z^2}(t,1) \right) dw\,ds \right) \quad (4.29)$$

Consequently, we obtain from (4.26) for $i = N-1$ and $i = 1$, (4.14), (4.18) and (4.28), (4.29) that for every $t \geq t_0$, $N > 1/h^*$, the following equalities hold:

$$x_1(t+\delta) = \left( 1 - 2\lambda a_1 + \lambda h^2 c_1 + \lambda \frac{4v_0}{3v_0 - 2hg_0} \left( a_1 - \frac{h}{2} b_1 \right) \right) x_1(t)$$
$$+ \lambda \left( a_1 + \frac{h}{2} b_1 - \frac{v_0}{3v_0 - 2hg_0} \left( a_1 - \frac{h}{2} b_1 \right) \right) x_2(t) + \lambda h^2 u_1(t) \quad (4.30)$$
$$- \lambda \left( a_1 - \frac{h}{2} b_1 \right) \frac{2hd_0(t)}{3v_0 - 2hg_0} + Q_0(\delta, h, t)$$

$$x_{N-1}(t+\delta) = \left( 1 - 2\lambda a_{N-1} + \lambda h^2 c_{N-1} + \lambda \left( a_{N-1} + \frac{h}{2} b_{N-1} \right) \frac{4v_1}{3v_1 + 2hg_1} \right) x_{N-1}(t)$$
$$+ \lambda \left( a_{N-1} - \frac{h}{2} b_{N-1} - \left( a_{N-1} + \frac{h}{2} b_{N-1} \right) \frac{v_1}{3v_1 + 2hg_1} \right) x_{N-2}(t) + \lambda h^2 u_{N-1}(t) \quad (4.31)$$
$$+ \lambda \left( a_{N-1} + \frac{h}{2} b_{N-1} \right) \frac{2hd_1(t)}{3v_1 + 2hg_1} + Q_1(\delta, h, t)$$

where

$$Q_0(\delta, h, t) := P(\delta, h, t, h)$$
$$+ \frac{\lambda v_0 (2a_1 - hb_1)}{2(3v_0 - 2hg_0)} \left( \int_0^{2h}\!\!\int_0^s \left( \frac{\partial^2 x}{\partial z^2}(t,w) - \frac{\partial^2 x}{\partial z^2}(t,0) \right) dw\,ds - 4 \int_0^h\!\!\int_0^s \left( \frac{\partial^2 x}{\partial z^2}(t,w) - \frac{\partial^2 x}{\partial z^2}(t,0) \right) dw\,ds \right) \quad (4.32)$$

$$Q_1(\delta, h, t) := P(\delta, h, t, 1-h)$$
$$+ \lambda \frac{v_1(2a_{N-1} + hb_{N-1})}{2(3v_1 + 2hg_1)} \left( \int_{1-2h}^{1}\!\!\int_s^1 \left( \frac{\partial^2 x}{\partial z^2}(t,w) - \frac{\partial^2 x}{\partial z^2}(t,1) \right) dw\,ds - 4 \int_{1-h}^{1}\!\!\int_s^1 \left( \frac{\partial^2 x}{\partial z^2}(t,w) - \frac{\partial^2 x}{\partial z^2}(t,1) \right) dw\,ds \right) \quad (4.33)$$

We next notice that the function $\omega : [0,1] \times [0,1] \to \Re$ given by (4.8), (4.9) is a continuous function, which by virtue of the fact that $a(z)\eta''(z) + b(z)\eta'(z) + c(z)\eta(z) = -\sigma\eta(z)$ for all $z \in [0,1]$ (recall Assumption (H4) and definitions (4.3)), satisfies

$$-2a_i + h^2 c_i + \left( a_i + \frac{h}{2} b_i \right) \frac{\eta_{i+1}}{\eta_i} + \left( a_i - \frac{h}{2} b_i \right) \frac{\eta_{i-1}}{\eta_i} = -\sigma h^2 (1 - \omega(ih, h)), \text{ for } h \in (0,1], \, i = 1,\ldots,N-1 \quad (4.34)$$

Let (arbitrary) $\varepsilon > 0$ be given and let $h^{**} \in (0,1)$ be selected in such a way that:

$$\left| \frac{\partial^2 x}{\partial z^2}(t,w) - \frac{\partial^2 x}{\partial z^2}(t,z) \right| \leq \varepsilon \text{ for all } t \in [t_0, T+1], w, z \in [0,1], |w-z| \leq 2h^{**} \quad (4.35)$$

$$\left| \frac{\partial x}{\partial t}(s,z) - \frac{\partial x}{\partial t}(t,z) \right| \leq \varepsilon, \text{ for all } z \in [0,1], s, t \in [t_0, T+1], |s-t| \leq h^{**} \quad (4.36)$$

$$|\omega(z,h)| \leq 1/2, \text{ for all } z \in [0,1] \text{ and } h \in (0, h^{**}) \quad (4.37)$$

where $\omega:[0,1]\times[0,1]\to\Re$ is the function defined by (4.8) and (4.9). The existence of $h^{**}\in(0,1)$ that satisfies (4.35), (4.36) follows from uniform continuity of the functions $\frac{\partial^2 x}{\partial z^2}(t,z)$ and $\frac{\partial x}{\partial t}(t,z)$ on the compact set $[t_0,T+1]\times[0,1]$. Indeed, notice that continuity of the mapping $z\to\frac{\partial^2 x}{\partial z^2}(t,z)$ for $(t,z)\in\Re_+\times[0,1]$ in conjunction with the continuity of the mappings $(t,z)\to\frac{\partial x}{\partial t}(t,z)$, $(t,z)\to x(t,z)$ for $(t,z)\in(0,+\infty)\times[0,1]$ and the fact that (2.6) holds for all $(t,z)\in(0,+\infty)\times(0,1)$, implies that (2.6) holds for all $(t,z)\in(0,+\infty)\times[0,1]$ and that the mapping $(t,z)\to\frac{\partial^2 x}{\partial z^2}(t,z)$ is continuous for $(t,z)\in(0,+\infty)\times[0,1]$. The existence of $h^{**}\in(0,1)$ that satisfies (4.37) follows from uniform continuity of the function $\omega:[0,1]\times[0,1]\to\Re$ defined by (4.8), (4.9) on the compact set $[0,1]\times[0,1]$. Let $N>1$ be an integer sufficiently large such that $N^{-1}<\min\left(h^*,h^{**},\frac{2}{3\sigma}\right)$ and define $h\in(0,1)$ by means of (4.18).

Next, let $m=\left[\frac{(1+2\bar{a}+C)(T-t_0)}{h^2}\right]+2$ and define $\lambda=\frac{T-t_0}{mh^2}$. Notice that (4.19) holds for this selection and that $T=t_0+m\delta$, where $\delta>0$ is defined by (4.18). Moreover, notice that the definitions of $h\in(0,1)$, $\delta>0$ in (4.18), (4.35), (4.36), (4.32), (4.33) and (4.25) guarantee that there exists a constant $S>0$ such that

$$\max_{i=1,\ldots n-1}\left(|P(\delta,h,t,ih)|\right)+|Q_0(\delta,h,t)|+|Q_1(\delta,h,t)|\le S\varepsilon\,\delta\,,\text{ for all }t\in[t_0,T+1] \tag{4.38}$$

Notice that by virtue of (4.14), the constant $S>0$ may be selected to be independent of $h\in(0,h^*)$, $\lambda\in(0,1)$, $\varepsilon>0$ and the time instants $0<t_0\le T$. Using (4.18), (4.26), (4.30), (4.31), (4.38), (4.14), (4.15), (4.16), the facts that $g_0\eta(0)+v_0\eta'(0)<0$ and $g_1\eta(1)+v_1\eta'(1)>0$ (recall assumption (H4)), the triangle inequality and (4.21), (4.22), (4.23), we get for all $t\in[t_0,T]$:

$$\begin{aligned}\frac{|x_1(t+\delta)|}{\eta_1}&\le\left(1-2\lambda a_1+\lambda h^2 c_1+\lambda\frac{4v_0}{3v_0-2hg_0}\left(a_1-\frac{h}{2}b_1\right)\right)\frac{|x_1(t)|}{\eta_1}\\&+\lambda\left(a_1+\frac{h}{2}b_1-\frac{v_0}{3v_0-2hg_0}\left(a_1-\frac{h}{2}b_1\right)\right)\frac{\eta_2}{\eta_1}\frac{|x_2(t)|}{\eta_2}\\&+\lambda\left(a_1-\frac{h}{2}b_1\right)\frac{2h|d_0(t)|}{\eta_1(3v_0-2hg_0)}+\delta\frac{S\varepsilon+|u_1(t)|}{\eta_1}\end{aligned} \tag{4.39}$$

$$\frac{|x_i(t+\delta)|}{\eta_i}\le\left(1-2\lambda a_i+\lambda h^2 c_i\right)\frac{|x_i(t)|}{\eta_i}+\lambda\left(a_i+\frac{h}{2}b_i\right)\frac{\eta_{i+1}}{\eta_i}\frac{|x_{i+1}(t)|}{\eta_{i+1}}+\lambda\left(a_i-\frac{h}{2}b_i\right)\frac{\eta_{i-1}}{\eta_i}\frac{|x_{i-1}(t)|}{\eta_{i-1}}+\delta\frac{S\varepsilon+|u_i(t)|}{\eta_i},$$
$$\text{for }i=2,\ldots,N-2 \tag{4.40}$$

$$\begin{aligned}\frac{|x_{N-1}(t+\delta)|}{\eta_{N-1}}&\le\left(1-2\lambda a_{N-1}+\lambda h^2 c_{N-1}+\lambda\left(a_{N-1}+\frac{h}{2}b_{N-1}\right)\frac{4v_1}{3v_1+2hg_1}\right)\frac{|x_{N-1}(t)|}{\eta_{N-1}}\\&+\lambda\left(a_{N-1}-\frac{h}{2}b_{N-1}-\left(a_{N-1}+\frac{h}{2}b_{N-1}\right)\frac{v_1}{3v_1+2hg_1}\right)\frac{\eta_{N-2}}{\eta_{N-1}}\frac{|x_{N-2}(t)|}{\eta_{N-2}}\\&+\lambda\left(a_{N-1}+\frac{h}{2}b_{N-1}\right)\frac{2h|d_1(t)|}{\eta_{N-1}(3v_1+2hg_1)}+\delta\frac{S\varepsilon+|u_{N-1}(t)|}{\eta_{N-1}}\end{aligned} \tag{4.41}$$

Consequently, using equalities (4.34) and definitions (4.18), (4.5) we obtain from (4.39), (4.40) and (4.41) for all $t\in[t_0,T]$:

$$\frac{|x_1(t+\delta)|}{\eta_1} \leq \left(1-\lambda\,\delta(1-\omega(h,h))\right)V_N(t) + \lambda\left(a_1 - \frac{h}{2}b_1\right)\frac{2h|d_0(t)|}{\eta_1(3v_0 - 2hg_0)} + \Gamma\varepsilon\,\delta$$
$$+ \lambda\left(a_1 - \frac{h}{2}b_1\right)\left(\frac{4v_0}{3v_0 - 2hg_0} - \frac{v_0}{3v_0 - 2hg_0}\frac{\eta_2}{\eta_1} - \frac{\eta_0}{\eta_1}\right)V_N(t) + \delta\|u\|$$
(4.42)

$$\frac{|x_i(t+\delta)|}{\eta_i} \leq \left(1-\sigma\,\delta(1-\omega(ih,h))\right)V_N(t) + (\|u\|+\Gamma\varepsilon)\,\delta, \quad i=2,\ldots,N-2$$
(4.43)

$$\frac{|x_{N-1}(t+\delta)|}{\eta_{N-1}} \leq \left(1-\sigma\,\delta(1-\omega(1-h,h))\right)V_N(t) + \lambda\left(a_{N-1} + \frac{h}{2}b_{N-1}\right)\frac{2h|d_1(t)|}{\eta_{N-1}(3v_1 + 2hg_1)} + \Gamma\varepsilon\,\delta$$
$$+ \lambda\left(a_{N-1} + \frac{h}{2}b_{N-1}\right)\left(\frac{4v_1}{3v_1 + 2hg_1} - \frac{v_1}{3v_1 + 2hg_1}\frac{\eta_{N-2}}{\eta_{N-1}} - \frac{\eta_N}{\eta_{N-1}}\right)V_N(t) + \delta\|u\|$$
(4.44)

where $\Gamma := \dfrac{S}{\min_{0\leq z\leq 1}(\eta(z))}$ and $\|u\| := \max_{0\leq s\leq T}\left(\max_{0\leq z\leq 1}\left(\dfrac{|u(s,z)|}{\eta(z)}\right)\right)$.

It is a matter of straightforward manipulations to show that estimates (4.42), (4.43), (4.44) in conjunction with definitions (4.18), (4.5) and inequalities (4.17) imply that the following implications hold for all $t\in[t_0, T]$:

(Q1) If $\max\left(\dfrac{\|d_0\|}{\gamma_0(h)}, \dfrac{\|d_1\|}{\gamma_1(h)}\right) \leq V_N(t)$ then $V_N(t+\delta) \leq \left(1-\sigma\,\delta\left(1-\max_{0\leq z\leq 1}(\omega(z,h))\right)\right)V_N(t) + (\|u\|+\Gamma\varepsilon)\,\delta$.

(Q2) If $\max\left(\dfrac{\|d_0\|}{\gamma_0(h)}, \dfrac{\|d_1\|}{\gamma_1(h)}\right) \geq V_N(t)$ then $V_N(t+\delta) \leq \left(1-\sigma\,\delta\left(1-\max_{0\leq z\leq 1}(\omega(z,h))\right)\right)\max\left(\dfrac{\|d_0\|}{\gamma_0(h)}, \dfrac{\|d_1\|}{\gamma_1(h)}\right) + (\|u\|+\Gamma\varepsilon)\,\delta$.

where $\|d_0\| := \max_{0\leq s\leq T}(|d_0(s)|)$, $\|d_1\| := \max_{0\leq s\leq T}(|d_1(s)|)$, $\gamma_0(h) := -g_0\eta(0) - v_0\eta'(0) + v_0\left(\eta'(0) - \dfrac{4\eta(h) - 3\eta(0) - \eta(2h)}{2h}\right)$,
$\gamma_1(h) := g_1\eta(1) + v_1\eta'(1) + v_1\left(\dfrac{3\eta(1) - 4\eta(1-h) + \eta(1-2h)}{2h} - \eta'(1)\right)$.

By virtue of (Q1) and (Q2) and the fact that $T = t_0 + m\delta$, it follows that properties (P1), (P2) of Lemma 4.1 hold for the sequences $y(j) := \max\left(\dfrac{\|d_0\|}{\gamma_0(h)}, \dfrac{\|d_1\|}{\gamma_1(h)}\right)$, $\varphi(j) = V_N(t_0 + j\delta)$ for $j=0,1,2,\ldots$ with $K := 1$, $a = \beta := 1 - \sigma\,\delta\left(1 - \max_{0\leq z\leq 1}(\omega(z,h))\right)$, $g := (\|u\|+\Gamma\varepsilon)\,\delta$. Notice that (4.37) in conjunction with (4.19) and the fact that $h = N^{-1} < \min\left(h^*, h^{**}, \dfrac{2}{3\sigma}\right)$ implies that $a \in (0,1)$. Therefore, we obtain from Lemma 4.1:

$$V_N(T) \leq \max\left(\max\left(\frac{\|d_0\|}{\gamma_0(h)}, \frac{\|d_1\|}{\gamma_1(h)}\right), \left(1 - \sigma\left(1 - \max_{0\leq z\leq 1}(\omega(z,h))\right)\delta\right)^m V_N(t_0)\right) + \frac{\Gamma\varepsilon + \|u\|}{\sigma\left(1 - \max_{0\leq z\leq 1}(\omega(z,h))\right)}$$
(4.45)

Using (4.45), the fact that $1 - \sigma\left(1 - \max_{0\leq z\leq 1}(\omega(z,h))\right)\delta \leq \exp\left(-\sigma\left(1 - \max_{0\leq z\leq 1}(\omega(z,h))\right)\delta\right)$ and the fact that $T = t_0 + m\delta$, we get (4.4). The proof is complete. ◁

We next proceed with the proof of Theorem 2.4.

**Proof of Theorem 2.4:** It suffices to prove the following estimate

$$\int_0^1 \sin(\pi z)|x(t,z)|dz \leq \exp\left(-\left(a\pi^2 - c\right)t\right) \int_0^1 \sin(\pi z)|x(0,z)|dz \qquad (4.46)$$

$$+ \frac{a\pi}{a\pi^2 - c}\left(\max_{0\leq s\leq t}(|d_0(s)|) + \exp\left(\frac{b}{2a}\right)\max_{0\leq s\leq t}(|d_1(s)|)\right) + \frac{1}{a\pi^2 - c}\max_{0\leq s\leq t}\left(\max_{0\leq z\leq 1}(|u(s,z)|\sin(\pi z))\right)$$

for the solution of

$$\frac{\partial x}{\partial t}(t,z) = a\frac{\partial^2 x}{\partial z^2}(t,z) + cx(t,z) + u(t,z) \text{ for } (t,z) \in (0,+\infty) \times (0,1) \qquad (4.47)$$

with

$$x(t,0) = d_0(t), \quad x(t,1) = \exp\left(\frac{b}{2a}\right)d_1(t) \text{ for } t \geq 0 \qquad (4.48)$$

where $c < a\pi^2$ is a constant. Notice that if $y(t,z)$ is a solution of (2.15), (2.16) then the function $x(t,z) = y(t,z)\exp\left(\frac{bz}{2a}\right)$ defined for $(t,z) \in \mathfrak{R}_+ \times [0,1]$ is a solution of (4.47), (4.48) with $c = k - \frac{b^2}{4a}$ and $u(t,z) := v(t,z)\exp\left(\frac{bz}{2a}\right)$.

In order to show estimate (4.46), it suffices to show that for every time instants $0 < t_0 \leq T$ and for every $\varepsilon > 0$, the following estimate holds for sufficiently large $N > 1$:

$$W_N(T) \leq \exp\left(1 - \left(a\pi^2(1 - r(N^{-1})) - c\right)(T - t_0)\right)W_N(t_0) + \frac{(1+a)\varepsilon + \|u\| + a\pi\left(\|d_0\| + \exp\left(\frac{b}{2a}\right)\|d_1\|\right)}{a\pi^2(1 - r(N^{-1})) - c}, \qquad (4.49)$$

where $\|d_0\| := \max_{0\leq s\leq T}(|d_0(s)|)$, $\|d_1\| := \max_{0\leq s\leq T}(|d_1(s)|)$, $\|u\| := \max_{0\leq s\leq T}\left(\max_{0\leq z\leq 1}(|u(s,z)|\sin(\pi z))\right)$,

$$r(h) := \begin{cases} 1 - 2\frac{1-\cos(\pi h)}{\pi^2 h^2} & \text{if } h \in (0,1] \\ 0 & \text{if } h = 0 \end{cases}, \qquad (4.50)$$

$$W_N(t) := N^{-1}\sum_{i=1}^{N-1}\sin(i\pi N^{-1})|x(t, iN^{-1})|. \qquad (4.51)$$

Notice that $r:[0,1] \to \mathfrak{R}$ as defined by (4.50) is a continuous function with $-1 \leq r(h) < 1$ for all $h \in [0,1]$. Indeed, (4.49) is sufficient for the proof of estimate (4.46), since

$$\int_{(i-1)h}^{ih}\sin(\pi z)|x(t,z)|dz \leq h\sin(i\pi h)|x_i(t)| + h^2 \max_{0\leq z\leq 1}\left(\left|\frac{\partial}{\partial z}(\sin(\pi z)x(t,z))\right|\right)$$

where $h := N^{-1}$, we obtain

$$\int_0^1 \sin(\pi z)|x(t,z)|dz \leq W_N(t) + h \max_{0\leq z\leq 1}\left(\left|\frac{\partial}{\partial z}(\sin(\pi z)x(t,z))\right|\right) \qquad (4.52)$$

Therefore, we obtain from (4.49) that

$$W_N(T) \leq \exp\left(1 - \left(a\pi^2(1 - r(N^{-1})) - c\right)(T - t_0)\right)W_N(t_0) + \frac{(1+a)\varepsilon + \|u\| + a\pi\left(\|d_0\| + \exp\left(\frac{b}{2a}\right)\|d_1\|\right)}{a\pi^2(1 - r(N^{-1})) - c} + 2MN^{-1}$$

where $M := \max_{0 \leq z \leq 1}\left(\left|\frac{\partial}{\partial z}(\sin(\pi z)x(T,z))\right|\right) + \max_{0 \leq z \leq 1}\left(\left|\frac{\partial}{\partial z}(\sin(\pi z)x(t_0,z))\right|\right)$. Letting $N \to +\infty$, we get (using (4.50)) for all $\varepsilon > 0$:

$$W_N(T) \leq \exp\left(1 - (a\pi^2 - c)(T - t_0)\right)W_N(t_0) + \frac{(1+a)\varepsilon + \|u\| + a\pi\left(\|d_0\| + \exp\left(\frac{b}{2a}\right)\|d_1\|\right)}{a\pi^2 - c}$$

Since the above inequality holds for all $t_0 \in (0,T]$, $\varepsilon > 0$ and since $\lim_{t_0 \to 0^+}\left(\int_0^1 \sin(\pi z)|x(t_0,z)|dz\right) = \int_0^1 \sin(\pi z)|x_0(z)|dz$, the above inequality implies that (4.46) holds for all $t > 0$. Inequality (4.46) holds for $t = 0$ as well.

Therefore, we have to prove that for every time instants $0 < t_0 \leq T$ and for every $\varepsilon > 0$, estimate (4.49) holds for sufficiently large $N > 1$.

Let $0 < t_0 \leq T$ be given (arbitrary) time instants. For every $h \in (0,1)$, $\delta > 0$, $t \geq t_0$, $z \in [h, 1-h]$ we obtain (4.11), (4.12). We next select an integer $N > 3$ and define:

$$h := N^{-1}, \ x_i(t) := x(t, ih), \ u_i(t) = u(t, ih), \ i = 0, \ldots, N, \ \delta := \lambda h^2 \tag{4.53}$$

where

$$\lambda \in \left(0, \frac{1}{1 + 2a + |c|}\right) \tag{4.54}$$

is a constant to be selected. Notice that (4.54) guarantee the following inequalities for all $N > 3$:

$$1 - 2\lambda a + \lambda h^2 c > 0 \tag{4.55}$$

It follows from (4.47), (4.11), (4.12) and the fact $\delta := \lambda h^2$ that for every $h \in (0,1)$, $\lambda \in \left(0, \frac{1}{1 + 2a + |c|}\right)$, $t \geq t_0$, $z \in [h, 1-h]$, the following equality holds:

$$x(t+\delta, z) = \left(1 - 2\lambda a + \lambda h^2 c\right)x(t,z) + \lambda a x(t, z+h) + \lambda a x(t, z-h) + \lambda h^2 u(t,z) + P(\delta, h, t, z) \tag{4.56}$$

where

$$P(\delta, h, t, z) := \int_t^{t+\delta}\left(\frac{\partial x}{\partial t}(s,z) - \frac{\partial x}{\partial t}(t,z)\right)ds - \lambda a \int_z^{z+h}\int_z^s\left(\frac{\partial^2 x}{\partial z^2}(t,w) - \frac{\partial^2 x}{\partial z^2}(t,z)\right)dw\,ds \\ -\lambda a \int_z^{z-h}\int_z^s\left(\frac{\partial^2 x}{\partial z^2}(t,w) - \frac{\partial^2 x}{\partial z^2}(t,z)\right)dw\,ds \tag{4.57}$$

Furthermore, it follows from (4.56) and definitions (4.53) that for every $t \geq t_0$, the following equalities hold:

$$x_i(t+\delta) = \left(1 - 2\lambda a + \lambda h^2 c\right)x_i(t) + \lambda a x_{i+1}(t) + \lambda a x_{i-1}(t) + \lambda h^2 u_i(t) + P(\delta, h, t, ih), \ i = 1, \ldots, N-1 \tag{4.58}$$

We next notice that the function $r:[0,1] \to \Re$ given by (4.50) is a continuous function, which satisfies

$$\frac{\sin(\pi i h - \pi h) + \sin(\pi i h + \pi h)}{\sin(\pi i h)} - 2 = -2(1 - \cos(\pi h)) = -\pi^2 h^2(1 - r(h)), \text{ for } i = 1, \ldots, N-1. \tag{4.59}$$

Let (arbitrary) $\varepsilon > 0$ be given and let $h^{**} \in (0,1)$ be selected in such a way that inequalities (4.35), (4.36) hold and

$$|r(h)| \leq 1/2, \text{ for all } h \in (0, h^{**}) \quad (4.60)$$

The existence of $h^{**} \in (0,1)$ that satisfies (4.35), (4.36) follows from uniform continuity of the functions $\frac{\partial^2 x}{\partial z^2}(t,z)$ and $\frac{\partial x}{\partial t}(t,z)$ on the compact set $[t_0, T+1] \times [0,1]$. Indeed, notice that continuity of the mapping $z \to \frac{\partial^2 x}{\partial z^2}(t,z)$ for $(t,z) \in \Re_+ \times [0,1]$ in conjunction with the continuity of the mappings $(t,z) \to \frac{\partial x}{\partial t}(t,z)$, $(t,z) \to x(t,z)$ for $(t,z) \in (0,+\infty) \times [0,1]$ and the fact that (4.47) holds for all $(t,z) \in (0,+\infty) \times (0,1)$, implies that (4.47) holds for all $(t,z) \in (0,+\infty) \times [0,1]$ and that the mapping $(t,z) \to \frac{\partial^2 x}{\partial z^2}(t,z)$ is continuous for $(t,z) \in (0,+\infty) \times [0,1]$. The existence of $h^{**} \in (0,1)$ that satisfies (4.60) follows from uniform continuity of the function $r : [0,1] \to \Re$ defined by (4.50) on the compact set $[0,1]$. Let $N > 1$ be an integer sufficiently large such that $N^{-1} < \min(h^{**}, 1/3)$ and define $h \in (0,1)$ by (4.53).

Next, let $m = \left[\frac{(1 + 2a + |c|)(T - t_0)}{h^2}\right] + 2$ and define $\lambda = \frac{T - t_0}{mh^2}$. Notice that (4.54) holds for this selection and that $T = t_0 + m\delta$, where $\delta > 0$ is defined by (4.53). Moreover, notice that the definitions of $h \in (0,1)$, $\delta > 0$ in (4.53), (4.35), (4.36) and (4.57) guarantee that there exists a constant $S > 0$ such that

$$\max_{i=1,\ldots,n-1} \left(|P(\delta, h, t, ih)|\right) \leq (1+a)\varepsilon \delta, \text{ for all } t \in [t_0, T+1] \quad (4.61)$$

Using (4.48), (4.53), (4.58), (4.61), the triangle inequality and (4.55), we get for all $t \in [t_0, T]$:

$$\sum_{i=1}^{N-1} \sin(i\pi h)|x_i(t+\delta)| \leq \sum_{i=1}^{N-1}\left(1 - 2\lambda a + \lambda h^2 c + \lambda a \frac{\sin((i+1)\pi h) + \sin((i-1)\pi h)}{\sin(i\pi h)}\right)\sin(i\pi h)|x_i(t)|$$
$$+ (1+a)\varepsilon \delta \sum_{i=1}^{N-1}\sin(i\pi h) + \lambda a \sin(\pi h)|d_0(t)| + \lambda a \sin((N-1)\pi h)\exp\left(\frac{b}{2a}\right)|d_1(t)| + \lambda h^2 \sum_{i=1}^{N-1}|u_i(t)|\sin(i\pi h) \quad (4.62)$$

Consequently, using definitions (4.53), (4.51), we obtain from (4.62) and (4.59) for all $t \in [t_0, T]$:

$$W_N(t+\delta) \leq \left(1 - \lambda a\pi^2 h^2 (1-r(h)) + \lambda h^2 c\right)W_N(t) + \lambda h^2 \max_{0 \leq z \leq 1}\left(|u(t,z)|\sin(\pi z)\right)$$
$$+ (1+a)\varepsilon \delta \sum_{i=1}^{N-1} h\sin(i\pi h) + \lambda ah\sin(\pi h)\left(\|d_0\| + \exp\left(\frac{b}{2a}\right)\|d_1\|\right) \quad (4.63)$$

where $\|d_0\| := \max_{0 \leq s \leq T}\left(|d_0(s)|\right)$, $\|d_1\| := \max_{0 \leq s \leq T}\left(|d_1(s)|\right)$. Applying (4.63) inductively for the sequence $W_N(t_0 + j\delta)$, for $j = 0, 1, \ldots, m$ and using the fact that $T = t_0 + m\delta$ and the facts $\sin(\pi h) \leq \pi h$, $\delta := \lambda h^2$, we get:

$$W_N(T) \leq \left(1 - \delta\left(a\pi^2(1-r(h)) - c\right)\right)^m W_N(t_0) + \frac{(1+a)\varepsilon + \|u\|}{a\pi^2(1-r(h)) - c} + \frac{a\pi}{a\pi^2(1-r(h)) - c}\left(\|d_0\| + \exp\left(\frac{b}{2a}\right)\|d_1\|\right)$$

where $\|u\| := \max_{0\le s\le T}\left(\max_{0\le z\le 1}(|u(s,z)||\sin(\pi z)|)\right)$. The above inequality in conjunction with the facts that $1-\delta(a\pi^2(1-r(h))-c) \le \exp(-\delta(a\pi^2(1-r(h))-c))$ and $T = t_0 + m\delta$ implies (4.49).

The proof is complete. ◁

Next, we give the proof of Theorem 2.3.

**Proof of Theorem 2.3:** Since the eigenfunctions $\{\phi_n\}_{n=1}^{\infty}$ of the Sturm-Liouville operator $A: D \to L_r^2([0,1])$ defined by (2.1), (2.2) under Assumptions (H1), (H2) form an orthonormal basis of $L_r^2([0,1])$, it follows that Parseval's identity holds, i.e.,

$$\|x[t]\|_{2,r}^2 = \sum_{n=1}^{\infty} c_n^2(t), \text{ for all } t \ge 0 \tag{4.64}$$

where

$$c_n(t) := \int_0^1 r(z) x(t,z) \phi_n(z) dz, \text{ for } n = 1, 2, \ldots \tag{4.65}$$

By virtue of (2.6), it follows from repeated integration by parts, that the following equalities hold for all $t > 0$:

$$\dot{c}_n(t) = \int_0^1 r(z) \frac{\partial x}{\partial t}(t,z) \phi_n(z) dz = \int_0^1 \frac{\partial}{\partial z}\left(p(z) \frac{\partial x}{\partial z}(t,z)\right) \phi_n(z) dz - \int_0^1 q(z) x(t,z) \phi_n(z) dz + \int_0^1 r(z) u(t,z) \phi_n(z) dz$$

$$= p(1) \frac{\partial x}{\partial z}(t,1) \phi_n(1) - p(0) \frac{\partial x}{\partial z}(t,0) \phi_n(0) - \int_0^1 p(z) \frac{\partial x}{\partial z}(t,z) \frac{d\phi_n}{dz}(z) dz$$

$$- \int_0^1 q(z) x(t,z) \phi_n(z) dz + \int_0^1 r(z) u(t,z) \phi_n(z) dz =$$

$$= p(1)\left(\frac{\partial x}{\partial z}(t,1)\phi_n(1) - x(t,1)\frac{d\phi_n}{dz}(1)\right) + p(0)\left(\frac{d\phi_n}{dz}(0) x(t,0) - \frac{\partial x}{\partial z}(t,0)\phi_n(0)\right)$$

$$+ \int_0^1 x(t,z)\left[\frac{d}{dz}\left(p(z)\frac{d\phi_n}{dz}(z)\right) - q(z)\phi_n(z)\right] dz + \int_0^1 r(z) u(t,z) \phi_n(z) dz =$$

$$= p(1)\left(\frac{\partial x}{\partial z}(t,1)\phi_n(1) - x(t,1)\frac{d\phi_n}{dz}(1)\right) + p(0)\left(\frac{d\phi_n}{dz}(0) x(t,0) - \frac{\partial x}{\partial z}(t,0)\phi_n(0)\right)$$

$$- \int_0^1 r(z) x(t,z) (A\phi_n)(z) dz + \int_0^1 r(z) u(t,z) \phi_n(z) dz$$

Thus we get for all $t > 0$:

$$\dot{c}_n(t) = p(1)\left(\frac{\partial x}{\partial z}(t,1)\phi_n(1) - x(t,1)\frac{d\phi_n}{dz}(1)\right) + p(0)\left(\frac{d\phi_n}{dz}(0) x(t,0) - \frac{\partial x}{\partial z}(t,0)\phi_n(0)\right)$$
$$- \int_0^1 r(z) x(t,z) (A\phi_n)(z) dz + \int_0^1 r(z) u(t,z) \phi_n(z) dz \tag{4.66}$$

It follows from (4.66), the fact that $(A\phi_n)(z) = \lambda_n \phi_n(z)$ and definition (4.65) that the following equation holds for all $t > 0$:

$$\dot{c}_n(t) + \lambda_n c_n(t) = p(1)\left(\frac{\partial x}{\partial z}(t,1)\phi_n(1) - x(t,1)\frac{d\phi_n}{dz}(1)\right)$$
$$+ p(0)\left(\frac{d\phi_n}{dz}(0) x(t,0) - \frac{\partial x}{\partial z}(t,0)\phi_n(0)\right) + \int_0^1 r(z) u(t,z) \phi_n(z) dz \tag{4.67}$$

Next, we show that for all $t \geq 0$:

$$x(t,0)\frac{d\phi_n}{dz}(0) - \frac{\partial x}{\partial z}(t,0)\phi_n(0) = \frac{d_0(t)}{g_0^2 + v_0^2}\left(g_0\frac{d\phi_n}{dz}(0) - v_0\phi_n(0)\right) \quad (4.68)$$

Indeed, the equation $g_0 x(t,0) + v_0 \frac{\partial x}{\partial z}(t,0) = d_0(t)$ (see (2.7)) gives:

$$g_0^2 x(t,0)\frac{d\phi_n}{dz}(0) + g_0 v_0 \frac{d\phi_n}{dz}(0)\frac{\partial x}{\partial z}(t,0) = d_0(t)g_0\frac{d\phi_n}{dz}(0)$$

$$g_0 v_0 \phi_n(0) x(t,0) + v_0^2 \phi_n(0)\frac{\partial x}{\partial z}(t,0) = d_0(t)v_0\phi_n(0)$$

from which we obtain:

$$d_0(t)\left(g_0\frac{d\phi_n}{dz}(0) - v_0\phi_n(0)\right) = (g_0^2 + v_0^2)\left(x(t,0)\frac{d\phi_n}{dz}(0) - \phi_n(0)\frac{\partial x}{\partial z}(t,0)\right)$$

$$+ \left(g_0\frac{\partial x}{\partial z}(t,0) - v_0 x(t,0)\right)\left(v_0\frac{d\phi_n}{dz}(0) + g_0\phi_n(0)\right)$$

Equation (4.68) follows directly from the above equation and the fact that $g_0\phi_n(0) + v_0\frac{d\phi_n}{dz}(0) = 0$.

Moreover, the equation $g_1 x(t,1) + v_1 \frac{\partial x}{\partial z}(t,1) = d_1(t)$ (see (2.7)) gives:

$$g_1^2 x(t,1)\frac{d\phi_n}{dz}(1) + g_1 v_1 \frac{d\phi_n}{dz}(1)\frac{\partial x}{\partial z}(t,1) = d_1(t)g_1\frac{d\phi_n}{dz}(1)$$

$$g_1 v_1 \phi_n(1) x(t,1) + v_1^2 \phi_n(1)\frac{\partial x}{\partial z}(t,1) = d_1(t)v_1\phi_n(1)$$

from which we obtain:

$$d_1(t)\left(g_1\frac{d\phi_n}{dz}(1) - v_1\phi_n(1)\right) = (g_1^2 + v_1^2)\left(x(t,1)\frac{d\phi_n}{dz}(1) - \phi_n(1)\frac{\partial x}{\partial z}(t,1)\right)$$

$$+ \left(g_1\frac{\partial x}{\partial z}(t,1) - v_1 x(t,1)\right)\left(v_1\frac{d\phi_n}{dz}(1) + g_1\phi_n(1)\right)$$

The fact that $g_1\phi_n(1) + v_1\frac{d\phi_n}{dz}(1) = 0$ in conjunction with the above equation implies that:

$$\phi_n(1)\frac{\partial x}{\partial z}(t,1) - x(t,1)\frac{d\phi_n}{dz}(1) = \frac{d_1(t)}{g_1^2 + v_1^2}\left(v_1\phi_n(1) - g_1\frac{d\phi_n}{dz}(1)\right) \quad (4.69)$$

Using (4.67), (4.68) and (4.69), we obtain for all $t > 0$:

$$\dot{c}_n(t) + \lambda_n c_n(t) = \frac{p(1)d_1(t)}{g_1^2 + v_1^2}\left(v_1\phi_n(1) - g_1\frac{d\phi_n}{dz}(1)\right)$$

$$+ \frac{p(0)d_0(t)}{g_0^2 + v_0^2}\left(g_0\frac{d\phi_n}{dz}(0) - v_0\phi_n(0)\right) + \int_0^1 r(z)u(t,z)\phi_n(z)dz \quad (4.70)$$

Integrating the differential equations (4.70), we obtain for all $0 < T \leq t$ and $n = 1,2,...$:

$$c_n(t) = \exp(-\lambda_n(t-T))c_n(T) + \frac{p(1)}{g_1^2 + v_1^2}\left(v_1\phi_n(1) - g_1\frac{d\phi_n}{dz}(1)\right)\int_T^t \exp(-\lambda_n(t-s))d_1(s)ds$$

$$+ \frac{p(0)}{g_0^2 + v_0^2}\left(g_0\frac{d\phi_n}{dz}(0) - v_0\phi_n(0)\right)\int_T^t \exp(-\lambda_n(t-s))d_0(s)ds + \int_T^t \exp(-\lambda_n(t-s))\left(\int_0^1 r(z)u(s,z)\phi_n(z)dz\right)ds \quad (4.71)$$

Continuity of the mapping $\Re_+ \ni T \to c_n(T)$ and (4.71) implies the following equations for all $t \geq 0$ and $n = 1,2,...$:

$$c_n(t) = \exp(-\lambda_n t)c_n(T) + \frac{p(1)}{g_1^2 + v_1^2}\left(v_1\phi_n(1) - g_1\frac{d\phi_n}{dz}(1)\right)\int_0^t \exp(-\lambda_n(t-s))d_1(s)ds$$

$$+ \frac{p(0)}{g_0^2 + v_0^2}\left(g_0\frac{d\phi_n}{dz}(0) - v_0\phi_n(0)\right)\int_0^t \exp(-\lambda_n(t-s))d_0(s)ds + \int_0^t \exp(-\lambda_n(t-s))\left(\int_0^1 r(z)u(s,z)\phi_n(z)dz\right)ds$$

(4.72)

Equations (4.72) in conjunction with the previous inequality imply the following estimates for all $t \geq 0$ and $n = 1,2,\ldots$:

$$|c_n(t)| \leq \exp(-\lambda_n t)|c_n(0)| + \frac{p(1)}{g_1^2 + v_1^2}\left|v_1\phi_n(1) - g_1\frac{d\phi_n}{dz}(1)\right|\frac{1-\exp(-\lambda_n t)}{\lambda_n}\max_{0\leq s\leq t}(|d_1(s)|)$$

$$+ \frac{p(0)}{g_0^2 + v_0^2}\left|g_0\frac{d\phi_n}{dz}(0) - v_0\phi_n(0)\right|\frac{1-\exp(-\lambda_n t)}{\lambda_n}\max_{0\leq s\leq t}(|d_0(s)|) + \frac{1-\exp(-\lambda_n t)}{\lambda_n}\|u\|\int_0^1 r(z)|\phi_n(z)|dz$$

(4.73)

where $\|u\| := \max_{0\leq s\leq t}\left(\max_{0\leq z\leq 1}(|u(s,z)|)\right)$. We next use the inequality $(a+b)^2 \leq (1+\varepsilon^{-1})a^2 + (1+\varepsilon)b^2$ (which holds for every $a,b \in \Re$ and $\varepsilon > 0$) in conjunction with (4.73). Using (4.73) for $t > 0$ with $\varepsilon = \frac{(1-\exp(-\lambda_n t))^2}{1-(1-\exp(-\lambda_n t))^2} > 0$, we get all $t > 0$ and $n = 1,2,\ldots$:

$$|c_n(t)|^2 \leq \frac{\exp(-\lambda_n t)}{2-\exp(-\lambda_n t)}|c_n(0)|^2$$

$$+ \frac{1}{\lambda_n^2}\left(\frac{p(1)}{g_1^2+v_1^2}\left|v_1\phi_n(1) - g_1\frac{d\phi_n}{dz}(1)\right|\max_{0\leq s\leq t}(|d_1(s)|) + \frac{p(0)}{g_0^2+v_0^2}\left|g_0\frac{d\phi_n}{dz}(0) - v_0\phi_n(0)\right|\max_{0\leq s\leq t}(|d_0(s)|) + \|u\|\int_0^1 r(z)|\phi_n(z)|dz\right)^2$$

(4.74)

Again, we next use the inequality $(a+b)^2 \leq (1+\varepsilon^{-1})a^2 + (1+\varepsilon)b^2$ (which holds for every $a,b \in \Re$ and $\varepsilon > 0$) in conjunction with (4.74). We obtain for all $t > 0$, $\varepsilon, \omega > 0$ and $n = 1,2,\ldots$:

$$|c_n(t)|^2 \leq \frac{\exp(-\lambda_n t)}{2-\exp(-\lambda_n t)}|c_n(0)|^2 + \frac{(1+\varepsilon^{-1})(1+\omega)}{\lambda_n^2}\left(\frac{p(0)}{g_0^2+v_0^2}\right)^2\left|g_0\frac{d\phi_n}{dz}(0) - v_0\phi_n(0)\right|^2\max_{0\leq s\leq t}(|d_0(s)|^2)$$

$$+ \frac{(1+\varepsilon)(1+\omega)}{\lambda_n^2}\left(\frac{p(1)}{g_1^2+v_1^2}\right)^2\left|v_1\phi_n(1) - g_1\frac{d\phi_n}{dz}(1)\right|^2\max_{0\leq s\leq t}(|d_1(s)|^2) + \frac{(1+\omega^{-1})}{\lambda_n^2}\|u\|^2\left(\int_0^1 r(z)|\phi_n(z)|dz\right)^2$$

(4.75)

We notice that (4.75) holds for $t = 0$ as well. Since $\lambda_n \geq \lambda_1 > 0$ for all $n = 1,2,\ldots$, we obtain from (4.75) the following estimates for all $t \geq 0$, $\varepsilon, \omega > 0$ and $n = 1,2,\ldots$:

$$|c_n(t)|^2 \leq \frac{\exp(-\lambda_1 t)}{2-\exp(-\lambda_1 t)}|c_n(0)|^2 + \frac{(1+\varepsilon^{-1})(1+\omega)}{\lambda_n^2}\left(\frac{p(0)}{g_0^2+v_0^2}\right)^2\left|g_0\frac{d\phi_n}{dz}(0) - v_0\phi_n(0)\right|^2\max_{0\leq s\leq t}(|d_0(s)|^2)$$

$$+ \frac{(1+\varepsilon)(1+\omega)}{\lambda_n^2}\left(\frac{p(1)}{g_1^2+v_1^2}\right)^2\left|v_1\phi_n(1) - g_1\frac{d\phi_n}{dz}(1)\right|^2\max_{0\leq s\leq t}(|d_1(s)|^2) + \frac{(1+\omega^{-1})}{\lambda_n^2}\|u\|^2\left(\int_0^1 r(z)|\phi_n(z)|dz\right)^2$$

(4.76)

Therefore, by virtue of estimates (4.76), (4.64), the following estimate holds for all $t \geq 0$ and $\varepsilon, \omega > 0$:

$$\|x[t]\|_r \leq \sqrt{\frac{\exp(-\lambda_1 t)}{2-\exp(-\lambda_1 t)}}\|x[0]\|_r + C_0\sqrt{(1+\varepsilon^{-1})(1+\omega)} \max_{0\leq s\leq t}(|d_0(s)|)$$
$$+ C_1\sqrt{(1+\varepsilon)(1+\omega)} \max_{0\leq s\leq t}(|d_1(s)|) + \tilde{C}\sqrt{1+\omega^{-1}} \max_{0\leq s\leq t}\left(\max_{0\leq z\leq 1}(|u(s,z)|)\right) \quad (4.77)$$

where

$$C_0 := \frac{p(0)}{\sqrt{g_0^2+v_0^2}} \sqrt{\sum_{n=1}^{\infty} \frac{1}{\lambda_n^2} \left|\frac{g_0}{\sqrt{g_0^2+v_0^2}} \frac{d\phi_n}{dz}(0) - \frac{v_0}{\sqrt{g_0^2+v_0^2}} \phi_n(0)\right|^2} \quad (4.78)$$

$$C_1 := \frac{p(1)}{\sqrt{g_1^2+v_1^2}} \sqrt{\sum_{n=1}^{\infty} \frac{1}{\lambda_n^2} \left|\frac{v_1}{\sqrt{g_1^2+v_1^2}} \phi_n(1) - \frac{g_1}{\sqrt{g_1^2+v_1^2}} \frac{d\phi_n}{dz}(1)\right|^2} \quad (4.79)$$

$$\tilde{C} := \sqrt{\sum_{n=1}^{\infty} \frac{1}{\lambda_n^2} \left(\int_0^1 r(z)|\phi_n(z)|dz\right)^2} \quad (4.80)$$

We notice that since $\int_0^1 r(z)\phi_n^2(z)dz = 1$ for $n=1,2,\ldots$, we get $\max_{0\leq z\leq 1}(\phi_n^2(z))\int_0^1 r(z)dz \geq 1$ and consequently $\max_{0\leq z\leq 1}(|\phi_n(z)|) \geq \left(\sqrt{\int_0^1 r(z)dz}\right)^{-1}$ for $n=1,2,\ldots$. Therefore, inequality (2.4) implies that $\sum_{n=1}^{\infty} \lambda_n^{-1} < +\infty$ and thus $\sum_{n=1}^{\infty} \lambda_n^{-2} < +\infty$. Using the Cauchy-Schwarz inequality and the fact that $\int_0^1 r(z)\phi_n^2(z)dz = 1$ for $n=1,2,\ldots$, we get $\left(\int_0^1 r(z)|\phi_n(z)|dz\right)^2 \leq \int_0^1 r(z)dz$ for $n=1,2,\ldots$. Thus we get:

$$\tilde{C} := \sqrt{\sum_{n=1}^{\infty} \frac{1}{\lambda_n^2} \left(\int_0^1 r(z)|\phi_n(z)|dz\right)^2} \leq \sqrt{\sum_{n=1}^{\infty} \lambda_n^{-2}} \sqrt{\int_0^1 r(z)dz} < +\infty$$

Using Lemma 2.1 in [14], it follows that the boundary value problem

$$\frac{d}{dz}\left(p(z)\frac{d\tilde{x}}{dz}(z)\right) - q(z)\tilde{x}(z) = 0, \text{ for all } z \in [0,1], \quad (4.81)$$

with

$$g_0\tilde{x}(0) + v_0\frac{d\tilde{x}}{dz}(0) = \sqrt{g_0^2+v_0^2} \quad , \quad g_1\tilde{x}(1) + v_1\frac{d\tilde{x}}{dz}(1) = 0 \quad (4.82)$$

has a unique solution $\tilde{x} \in C^2([0,1])$, which satisfies

$$p^2(0)\sum_{n=1}^{\infty} \lambda_n^{-2} \left|\frac{g_0}{\sqrt{g_0^2+v_0^2}} \frac{d\phi_n}{dz}(0) - \frac{v_0}{\sqrt{g_0^2+v_0^2}} \phi_n(0)\right|^2 = \int_0^1 r(z)\tilde{x}^2(z)dz \quad (4.83)$$

Changing $z$ by $1-z$ and using Lemma 2.1 in [14], it follows that the boundary value problem

$$\frac{d}{dz}\left(p(z)\frac{d\bar{x}}{dz}(z)\right) - q(z)\bar{x}(z) = 0, \text{ for all } z \in [0,1], \quad (4.84)$$

with

$$g_0 \bar{x}(0) + v_0 \frac{d\bar{x}}{dz}(0) = 0 \quad , \quad g_1 \bar{x}(1) + v_1 \frac{d\bar{x}}{dz}(1) = \sqrt{g_1^2 + v_1^2} \qquad (4.85)$$

has a unique solution $\bar{x} \in C^2([0,1])$, which satisfies

$$p^2(1) \sum_{n=1}^{\infty} \lambda_n^{-2} \left| \frac{g_1}{\sqrt{g_1^2 + v_1^2}} \frac{d\phi_n}{dz}(1) - \frac{v_1}{\sqrt{g_1^2 + v_1^2}} \phi_n(1) \right|^2 = \int_0^1 r(z) \bar{x}^2(z) dz \qquad (4.86)$$

The proof is complete. ◁

**Proof of Corollary 2.5:** Inequalities (2.20) and (2.22) are straightforward applications of Theorem 2.2 (with $\sigma \in (0, a\pi^2)$, $\theta \in (0, \pi - \varphi)$, $\varphi := \sqrt{\frac{\sigma}{a}}$ and $\eta(z) = \sin(\theta + z\varphi)$) and Theorem 2.4. Inequality (2.21) is a consequence of a similar analysis to that presented in the proof of Theorem 2.3, which is presented next.

Since the functions $\{\sqrt{2}\sin(n\pi z)\}_{n=1}^{\infty}$ form an orthonormal basis of $L^2([0,1])$, it follows that Parseval's identity holds, i.e.,

$$\int_0^1 x^2(t,z) dz = \sum_{n=1}^{\infty} c_n^2(t), \text{ for all } t \geq 0 \qquad (4.87)$$

where

$$c_n(t) := \sqrt{2} \int_0^1 x(t,z) \sin(n\pi z) dz, \text{ for } n = 1,2,... \qquad (4.88)$$

By virtue of (2.16), (2.19), it follows from repeated integration by parts, that the following equalities hold for all $t > 0$:

$$\dot{c}_n(t) = -a(-1)^n n\pi \sqrt{2} d_1(t) + an\pi \sqrt{2} d_0(t) - an^2 \pi^2 c_n(t)$$

Integrating the above differential equations, we obtain for all $0 < T \leq t$ and $n = 1,2,...$:

$$c_n(t) = \exp(-an^2 \pi^2 (t-T)) c_n(T) + an\pi \sqrt{2} \int_T^t \exp(-an^2 \pi^2 (t-s)) (d_0(s) - (-1)^n d_1(s)) ds \qquad (4.89)$$

Continuity of the mapping $\Re_+ \ni T \to c_n(T)$ and (4.89) implies the following equations for all $t \geq 0$ and $n = 1,2,...$:

$$c_n(t) = \exp(-an^2 \pi^2 t) c_n(0) + an\pi \sqrt{2} \int_0^t \exp(-an^2 \pi^2 (t-s)) (d_0(s) - (-1)^n d_1(s)) ds \qquad (4.90)$$

Equations (4.90) imply the following estimates for all $t \geq 0$ and $n = 1,2,...$:

$$|c_n(t)| \leq \exp(-an^2 \pi^2 t) c_n(0) + \sqrt{2} \frac{1 - \exp(-an^2 \pi^2 t)}{n\pi} \left( \max_{0 \leq s \leq t} (|d_0(s)|) + \max_{0 \leq s \leq t} (|d_1(s)|) \right) \qquad (4.91)$$

Therefore, we obtain from (4.91) the following estimates for all $t \geq 0$, $\varepsilon_n(t) > 0$ and $n = 1,2,...$:

$$|c_n(t)|^2 \leq (1 + \varepsilon_n^{-1}(t)) \exp(-2an^2 \pi^2 t) |c_n(0)|^2$$
$$+ 2(1 + \varepsilon_n(t)) \frac{(1 - \exp(-an^2 \pi^2 t))^2}{n^2 \pi^2} \left( \max_{0 \leq s \leq t} (|d_1(s)|) + \max_{0 \leq s \leq t} (|d_0(s)|) \right)^2 \qquad (4.92)$$

Selecting $\varepsilon_n(t) = \frac{1 - (1 - \exp(-an^2 \pi^2 t))^2}{(1 - \exp(-an^2 \pi^2 t))^2}$ for $t > 0$, we obtain from (4.92) for all $t > 0$ and $n = 1,2,...$:

$$|c_n(t)|^2 \leq \frac{\exp(-an^2\pi^2 t)}{2-\exp(-an^2\pi^2 t)}|c_n(0)|^2 + \frac{2}{n^2\pi^2}\left(\max_{0\leq s\leq t}(|d_1(s)|) + \max_{0\leq s\leq t}(|d_0(s)|)\right)^2 \quad (4.93)$$

Using the inequality $\frac{\exp(-an^2\pi^2 t)}{2-\exp(-an^2\pi^2 t)} \leq \frac{\exp(-a\pi^2 t)}{2-\exp(-a\pi^2 t)}$ for all $t>0$, $n=1,2,...$ and the fact that

$\sum_{n=1}^{\infty}\frac{2}{n^2\pi^2} = \frac{1}{3}$ (since $\sqrt{2}\int_0^1 z\sin(n\pi z)dz = -\frac{\sqrt{2}}{n\pi}(-1)^n$ and $\int_0^1 z^2 dz = \frac{1}{3}$), we get (2.21) for $t>0$ by combining (4.93) with (4.87). Estimate (2.21) holds for $t=0$ as well.

The proof is complete. ◁

For the proofs of Theorem 3.2 and Theorem 3.3, we need the following technical lemma. Its proof is provided in the Appendix.

**Lemma 4.2:** *For every $\sigma>0$, $M\geq 1$, $\varepsilon>0$ there exist a constant $\delta\in(0,\sigma)$ with the following property: if $\varphi:\Re_+\to\Re_+$ and $y:\Re_+\to\Re_+$ are locally bounded functions for which there exists a constant $\gamma\geq 0$ such that the following inequality holds for all $t_0\geq 0$ and $t\geq t_0$*

$$\varphi(t) \leq M\exp(-\sigma(t-t_0))\varphi(t_0) + \gamma \sup_{t_0\leq s\leq t}(y(s)), \quad (4.94)$$

*then the following inequality holds for all $t\geq 0$:*

$$\varphi(t) \leq M\exp(-\delta t)\varphi(0) + \gamma(1+\varepsilon)\sup_{0\leq s\leq t}(y(s)\exp(-\delta(t-s))) \quad (4.95)$$

Since the proofs of Theorem 3.2 and Theorem 3.3 are similar, they are presented together.

**Proofs of Theorem 3.2 and Theorem 3.3:** Let $x_0 \in \Delta$ (arbitrary) and consider the unique solution $x \in C^0(\Re_+\times[0,1])\cap C^1((0,+\infty)\times[0,1])$ of the evolution equation (2.19) with (3.1) with $x[t]\in C^2([0,1])$ for all $t\geq 0$ and $x(0,z)=x_0(z)$ for all $z\in[0,1]$. By virtue of Corollary 2.5 the solution satisfies estimates (2.20), (2.21) for all $t\geq 0$ with

$$d_0(t):=\int_0^1 g_0(z)x(t,z)dz \text{ and } d_1(t)=\int_0^1 g_1(z)x(t,z)dz \text{ for all } t\geq 0 \quad (4.96)$$

Using Hölder's inequality and definitions (4.96), we obtain for all $t\geq 0$:

$$|d_0(t)|\leq \|x[t]\|_2 \sqrt{\int_0^1 g_0^2(z)dz} \qquad |d_0(t)|\leq \|x[t]\|_{1,w}\sup_{0<z<1}\left(\frac{|g_0(z)|}{\sin(\pi z)}\right)$$
$$\text{and} \qquad (4.97)$$
$$|d_1(t)|\leq \|x[t]\|_2 \sqrt{\int_0^1 g_1^2(z)dz} \qquad |d_1(t)|\leq \|x[t]\|_{1,w}\sup_{0<z<1}\left(\frac{|g_1(z)|}{\sin(\pi z)}\right)$$

Consequently, we get from (2.21) and (4.97) for all $t\geq 0$:

$$\|x[t]\|_2 \leq \exp\left(-\frac{a\pi^2}{2}t\right)\|x[0]\|_2 + \gamma_2 \max_{0\leq s\leq t}(\|x[s]\|_2), \quad (4.98)$$

$$\|x[t]\|_{1,w} \leq \exp\left(-\frac{a\pi^2}{2}t\right)\|x[0]\|_{1,w} + \gamma_{1,w} \max_{0\leq s\leq t}\left(\|x[s]\|_{1,w}\right), \tag{4.99}$$

where $\gamma_2 := \frac{1}{\sqrt{3}}\left(\sqrt{\int_0^1 g_0^2(z)dz} + \sqrt{\int_0^1 g_1^2(z)dz}\right)$ and $\gamma_{1,w} = \frac{1}{\pi}\left(\sup_{0<z<1}\left(\frac{|g_0(z)|}{\sin(\pi z)}\right) + \sup_{0<z<1}\left(\frac{|g_1(z)|}{\sin(\pi z)}\right)\right)$. By virtue of (3.4) or (3.6), we have that $\gamma < 1$, where $\gamma = \gamma_2$ in case of Theorem 3.2 and $\gamma = \gamma_{1,w}$ in case of Theorem 3.3. Consequently, there exists $\varepsilon > 0$ such that $\gamma(1+\varepsilon) < 1$. By virtue of Lemma 4.2, for the selected $\varepsilon > 0$ for which $\gamma(1+\varepsilon) < 1$, there exists a constant $\delta \in \left(0, \frac{a\pi^2}{2}\right]$ with the following property: if $\varphi : \Re_+ \to \Re_+$ and $y : \Re_+ \to \Re_+$ are locally bounded functions for which inequality (4.94) holds for all $t_0 \geq 0$ and $t \geq t_0$ with $\sigma = \frac{a\pi^2}{2}$ and $M = 1$, then inequality (4.95) holds. Applying this property to the functions $\varphi(t) = y(t) = \|x[t]\|_2$ in case of Theorem 3.2 and the functions $\varphi(t) = y(t) = \|x[t]\|_{1,w}$ in case of Theorem 3.3 and noticing that inequality (4.94) holds for these functions (by virtue of estimates (4.98), (4.99) and the semigroup property for the solution of the evolution equation (2.19) with (3.1)), we get for Theorem 3.2

$$\|x[t]\|_2 \exp(\delta t) \leq \|x[0]\|_2 + \gamma(1+\varepsilon)\max_{0\leq s\leq t}\left(\|x[s]\|_2 \exp(\delta s)\right), \text{ for all } t \geq 0 \tag{4.100}$$

and for Theorem 3.3

$$\|x[t]\|_{1,w} \exp(\delta t) \leq \|x[0]\|_{1,w} + \gamma(1+\varepsilon)\max_{0\leq s\leq t}\left(\|x[s]\|_{1,w} \exp(\delta s)\right), \text{ for all } t \geq 0 \tag{4.101}$$

Since $\gamma(1+\varepsilon) < 1$, estimate (3.5) (in case of Theorem 3.2) or (3.7) (in case of Theorem 3.3) follows from estimate (4.100) or estimate (4.101), respectively, with $M = \frac{1}{1-\gamma(1+\varepsilon)}$.

The proof is complete. ◁

For the proof of Theorem 3.1, we need the following technical lemma. Its proof is provided in the Appendix.

**Lemma 4.3:** *For every $\sigma > 0$, $M \geq 1$, $\varepsilon > 0$ there exist a constant $\delta \in (0, \sigma]$ with the following property: if $q : \Re_+ \to \Re_+$ and $y : \Re_+ \to \Re_+$ are locally bounded functions for which there exists a constant $\gamma \geq 0$ such that the following inequality holds for all $t_0 \geq 0$ and $t \geq t_0$*

$$q(t) \leq \max\left(M\exp(-\sigma(t-t_0))q(t_0), \gamma \sup_{t_0\leq s\leq t}(y(s))\right), \tag{4.102}$$

*then the following inequality holds for all $t \geq 0$:*

$$q(t) \leq \max\left(M\exp(-\delta t)q(0), \gamma(1+\varepsilon)\sup_{0\leq s\leq t}\left(y(s)\exp(-\delta(t-s))\right)\right) \tag{4.103}$$

Finally, we present the proof of Theorem 3.1.

**Proof of Theorem 3.1:** Let $x_0 \in \Delta$ (arbitrary) and consider the unique solution $x \in C^0(\Re_+ \times [0,1]) \cap C^1((0,+\infty) \times [0,1])$ of the evolution equation (2.19) with (3.1) with $x[t] \in C^2([0,1])$ for all $t \geq 0$ and $x(0,z) = x_0(z)$ for all $z \in [0,1]$. By virtue of Corollary 2.5 the solution satisfies estimates (2.22) for all $t \geq 0$ with $\varphi \in (0,\pi)$, $\theta \in (0, \pi - \varphi)$ being the constants involved in (3.2), $\sigma = a\varphi^2$ and $(d_0, d_1)$ given by (4.96). Using Hölder's inequality and definitions (4.96), we obtain for all $t \geq 0$:

$$|d_0(t)| \leq \|x[t]\|_{\infty,w} \int_0^1 |g_0(z)| \sin(\theta + z\varphi) dz$$

$$|d_1(t)| \leq \|x[t]\|_{\infty,w} \int_0^1 |g_1(z)| \sin(\theta + z\varphi) dz \qquad (4.104)$$

where $\|x[t]\|_{\infty,w} := \max_{0 \leq z \leq 1} \left( \frac{|x(t,z)|}{\sin(\theta + z\varphi)} \right)$. Consequently, we get from (2.22) and (4.104) for all $t \geq 0$:

$$\|x[t]\|_{\infty,w} \leq \max\left( \exp(-\sigma t) \|x[0]\|_{\infty,w}, \gamma \max_{0 \leq s \leq t} (\|x[s]\|_{\infty,w}) \right), \qquad (4.105)$$

where $\gamma = \max\left( \frac{1}{\sin(\theta)} \int_0^1 |g_0(z)| \sin(\theta + z\varphi) dz, \frac{1}{\sin(\theta + \varphi)} \int_0^1 |g_1(z)| \sin(\theta + z\varphi) dz \right)$. By virtue of (3.2) we have that $\gamma < 1$. Consequently, there exists $\varepsilon > 0$ such that $\gamma(1+\varepsilon) < 1$. By virtue of Lemma 4.3, for the selected $\varepsilon > 0$ for which $\gamma(1+\varepsilon) < 1$, there exists a constant $\delta \in (0, \sigma]$ with the following property: if $q : \Re_+ \to \Re_+$ and $y : \Re_+ \to \Re_+$ are locally bounded functions for which inequality (4.102) holds for all $t_0 \geq 0$ and $t \geq t_0$ with $M = 1$, then inequality (4.103) holds. Applying this property to the functions $q(t) = y(t) = \|x[t]\|_{\infty,w}$ and noticing that inequality (4.102) holds for these functions (by virtue of estimate (4.105) and the semigroup property for the solution of the evolution equation (2.19) with (3.1)), we get

$$\|x[t]\|_{\infty,w} \exp(\delta t) \leq \max\left( \|x[0]\|_{\infty,w}, \gamma(1+\varepsilon) \max_{0 \leq s \leq t} (\|x[s]\|_{\infty,w} \exp(\delta s)) \right), \text{ for all } t \geq 0 \qquad (4.106)$$

Since $\gamma(1+\varepsilon) < 1$ and since there exist constants $K_1, K_2 > 0$ such that $K_1 \|x[t]\|_\infty \leq \|x[t]\|_{\infty,w} \leq K_2 \|x[t]\|_\infty$, estimate (3.3) follows from estimate (4.106). The proof is complete. ◁

## 5. Concluding Remarks

We have studied 1-D parabolic PDEs with disturbances at both boundaries and distributed disturbances. Our main results provided ISS estimates in weighted $L^\infty$, $L^2$ and $L^1$ norms. Due to the lack of an ISS Lyapunov functional for boundary disturbances, we have utilized a finite-difference scheme for the proof of the ISS estimates in weighted $L^\infty$ and $L^1$ norms. The ISS estimate for the sup-norm led to a refinement of the well-known maximum principle for the heat equation. We also applied the obtained results to quasi-static thermoelasticity models that involve nonlocal boundary conditions. By means of small-gain arguments, we were able to derive conditions that guarantee the global exponential stability of the zero solution for such models.

Future work may involve the study of the robustness of the closed-loop system for 1-D parabolic PDEs under boundary feedback control. It should be noticed that the application boundary feedback leads to a closed-loop system which is described by nonlocal boundary conditions (like the systems that were studied in Section 3 of the present work; see [26]). Moreover, future work may also involve the derivation of ISS estimates in general $L^p$ norms, where $1 \leq p \leq \infty$.

**References**

[1] Argomedo, F. B., E. Witrant, and C. Prieur, "$D^1$-Input-to-State Stability of a Time-Varying Nonhomogeneous Diffusive Equation Subject to Boundary Disturbances", *Proceedings of the American Control Conference*, Montreal, QC, 2012, 2978 - 2983.

[2] Argomedo, F. B., C. Prieur, E. Witrant, and S. Bremond, "A Strict Control Lyapunov Function for a Diffusion Equation With Time-Varying Distributed Coefficients", *IEEE Transactions on Automatic Control,* 58(2.2), 2013, 290-303.

[3] Boyce, W. E., and R. C. Diprima, *Elementary Differential Equations and Boundary Value Problems*, 6th Edition, Wiley, 1997.

[4] Dashkovskiy, S., and A. Mironchenko, "On the Uniform Input-to-State Stability of Reaction-Diffusion Systems", *Proceedings of the 49th Conference on Decision and Control*, Atlanta, GA, USA, 2010, 6547-6552.

[5] Dashkovskiy, S., and A. Mironchenko, "Local ISS of Reaction-Diffusion Systems", *Proceedings of the 18th IFAC World Congress*, Milano, Italy, 2011, 11018-11023.

[6] Dashkovskiy, S., and A. Mironchenko, "Input-to-State Stability of Infinite-Dimensional Control Systems", *Mathematics of Control, Signals, and Systems*, 25(2.1), 2013, 1-35.

[7] Day, W. A., "Extension of a Property of the Heat Equation to Linear Thermoelasticity and Other Theories", *Quarterly of Applied Mathematics*, 40, 1982, 319-330.

[8] Day, W. A., "A Decreasing Property of Solutions of Parabolic Equations with Applications to Thermoelasticity", *Quarterly of Applied Mathematics*, 40, 1983, 468-475.

[9] Ekolin, G., "Finite Difference Methods for a Nonlocal Boundary Value Problem for the Heat Equation", *BIT*, 31, 1991, 245-261.

[10] Fairweather, G., and J. C. Lopez-Marcos, "Galerkin Methods for a Semilinear Parabolic Problem with Nonlocal Boundary Conditions", *Advances in Computational Mathematics*, 6(2.1), 1996, 243-262.

[11] John, F., *Partial Differential Equations*, 4th Ed., Springer-Verlag, New York, 1982.

[12] Karafyllis, I. and Z.-P. Jiang, *Stability and Stabilization of Nonlinear Systems*, Springer-Verlag London (Series: Communications and Control Engineering), 2011.

[13] Karafyllis, I., and M. Krstic, "On the Relation of Delay Equations to First-Order Hyperbolic Partial Differential Equations", *ESAIM Control, Optimisation and Calculus of Variations*, 20(2.3), 2014, 894 - 923.

[14] Karafyllis, I., and M. Krstic, "ISS With Respect to Boundary Disturbances for 1-D Parabolic PDEs", to appear in *IEEE Transactions on Automatic Control*.

[15] Liu, Y., "Numerical Solution of the Heat equation with Nonlocal Boundary Conditions", *Journal of Computational and Applied Mathematics*, 110, 1999, 115-127.

[16] Mazenc, F., and C. Prieur, "Strict Lyapunov Functions for Semilinear Parabolic Partial Differential Equations", *Mathematical Control and Related Fields*, AIMS, 1(2.2), 2011, 231-250.

[17] Mironchenko, A. and H. Ito, "Integral Input-to-State Stability of Bilinear Infinite-Dimensional Systems", *Proceedings of the 53rd IEEE Conference on Decision and Control*, Los Angeles, California, USA, 2014, 3155-3160.

[18] Mironchenko, A. and H. Ito, "Construction of Lyapunov Functions for Interconnected Parabolic Systems: An iISS Approach", arXiv:1410.3058 [math.DS].

[19] Mironchenko, A., "Local Input-to-State Stability: Characterizations and Counterexamples", *Systems and Control Letters*, 87, 2016, 23-28.

[20] Mironchenko, A. and F. Wirth, "Restatements of Input-to-State Stability in Infinite Dimensions: What Goes Wrong", submitted to the 22th International Symposium on Mathematical Theory of Systems and Networks (MTNS 2016).

[21] Mironchenko, A. and F. Wirth, "Global Converse Lyapunov Theorems for Infinite-Dimensional Systems", submitted to the *10th IFAC Symposium on Nonlinear Control Systems* (NOLCOS 2016).

[22] Naylor, A. W., and G. R. Sell, *Linear Operator Theory in Engineering and Science*, Springer, 1982.

[23] Pazy, A., *Semigroups of Linear Operators and Applications to Partial Differential Equations*, New York, Springer, 1983.

[24] Prieur, C., and F. Mazenc, "ISS-Lyapunov Functions for Time-Varying Hyperbolic Systems of Balance Laws", *Mathematics of Control, Signals, and Systems*, 24(1-2), 2012, 111-134.

[25] Smoller, J., *Shock Waves and Reaction-Diffusion Equations*, 2$^{nd}$ Edition, Springer-Verlag, New York, 1994.

[26] Smyshlyaev, A., and M. Krstic, "Closed-Form Boundary State Feedbacks for a Class of 1-D Partial Integro-Differential Equations", *IEEE Transactions on Automatic Control*, 49(2.13), 2004, 2185-2202.

[27] Sontag, E.D., "Smooth Stabilization Implies Coprime Factorization", *IEEE Transactions on Automatic Control*, 34, 1989, 435-443.

# Appendix

**Proof of Lemma 4.1:** Let arbitrary $j \in \{0,1,...,m\}$ be given. Define the set

$$A = \left\{ s \in \{0,...,j\} : y(s) > K\varphi(s) \right\} \tag{A.1}$$

Clearly, if $j \in A$ then we get $\varphi(j) \leq K^{-1} y(j) \leq K^{-1} \max_{s=0,...,j} (y(s))$ and estimate (4.2) holds in this case. Therefore, we next assume that $j \notin A$. We distinguish the following cases.

Case 1: $A = \varnothing$.
In this case we have $K\varphi(s) \leq y(s)$ for all $s = 0,...,j$. Property (P1) implies that $\varphi(s+1) \leq a\varphi(s) + g$ for all $s = 0,...,j$. Using induction, we get $\varphi(s) \leq a^s \varphi(0) + g \frac{1-a^s}{1-a}$ for all $s = 0,...,j$, which directly shows that estimate (4.2) holds in this case.

Case 2: $A \neq \varnothing$
Define $T = \max \{ s \in A \}$ and notice that since $j \notin A$, we obtain $T \leq j-1$. Definition (A.1) implies that $K\varphi(T) \leq y(T)$ and that $K\varphi(s) \geq y(s)$ for all $s = T+1,...,j$. Property (P1) implies that $\varphi(s+1) \leq a\varphi(s) + g$ for all $s = T+1,...,j$. Using induction, we get $\varphi(s) \leq a^{s-T-1} \varphi(T+1) + g \frac{1-a^{s-T-1}}{1-a}$ for all $s = T+1,...,j$. Moreover, Property (P2) implies that $\varphi(T+1) \leq \beta y(T) + g$, which combined with the previous inequality, implies that $\varphi(j) \leq a^{j-T-1} \beta y(T) + g \frac{1-a^{j-T}}{1-a}$. The last inequality shows that estimate (4.2) holds in this case as well. The proof is complete. ◁

**Proof of Lemma 4.2:** Let $\varepsilon > 0$ (arbitrary) be given and let constants $\lambda \in (0,1)$, $\mu > 1$ so that

$$0 < \frac{\mu}{1-\lambda\mu} \leq 1 + \varepsilon \tag{A.2}$$

Such constants $\lambda \in (0,1)$, $\mu > 1$ exist by virtue of continuity of the function $f(\lambda,\mu) := \frac{\mu}{1-\lambda\mu}$ in a neighborhood of $(0,1)$ and the fact that $f(0,1) = 1 < 1 + \varepsilon$. Define

$$T = \sigma^{-1} \ln(M/\lambda), \quad \omega := \sigma - T^{-1} \ln(M), \quad \delta := T^{-1} \ln(\mu) \tag{A.3}$$

Notice that definitions (A.3) guarantee that $\lambda = M \exp(-\sigma T) = \exp(-\omega T) < 1$ and that $\omega \in (0, \sigma]$, $\delta \in (0, \omega)$. Using (4.94), we obtain for all $p \geq 0$:

$$\varphi((i+1)T + p) \leq \lambda \varphi(iT + p) + \gamma \sup_{p+iT \leq s \leq p+(i+1)T} (y(s)), \text{ for all } i \geq 0 \tag{A.4}$$

Applying (A.4) repeatedly, we get:

$$\varphi(iT + p) \leq \lambda^i \varphi(p) + \gamma \sum_{k=0}^{i-1} \lambda^{i-1-k} \sup_{p+kT \leq s \leq p+(k+1)T} (y(s)), \text{ for all } i \geq 1 \text{ and } p \geq 0 \tag{A.5}$$

Using the facts that $\lambda = \exp(-\omega T)$ and $\delta \in (0, \omega)$, we get for all $i \geq 1$ and $p \geq 0$:

$$\varphi(iT + p) \leq \exp(-\omega iT) \varphi(p) + \gamma \exp(\omega T) \exp(-\delta iT) \sum_{k=0}^{i-1} \exp(-(\omega - \delta)(iT - kT)) \sup_{p+kT \leq s \leq p+(k+1)T} (y(s) \exp(\delta kT))$$

$$\leq \exp(-\omega iT) \varphi(p) + \gamma \exp(\omega T) \sum_{k=0}^{i-1} \exp(-(\omega - \delta)T(i-k)) \sup_{p+kT \leq s \leq p+(k+1)T} (y(s) \exp(-\delta(iT + p - s)))$$

$$\leq \exp(-\omega iT) \varphi(p) + \gamma \exp(\omega T) \sup_{p \leq s \leq p+iT} (y(s) \exp(-\delta(iT + p - s))) \sum_{k=0}^{i-1} \exp(-(\omega - \delta)T(i-k))$$

$$\leq \exp(-\omega iT) \varphi(p) + \gamma \exp(\omega T) \sup_{p \leq s \leq p+iT} (y(s) \exp(-\delta(iT + p - s))) \sum_{k=1}^{i} \exp(-(\omega - \delta)Tk)$$

$$\leq \exp(-\omega iT) \varphi(p) + \gamma \exp(\delta T) \sup_{p \leq s \leq p+iT} (y(s) \exp(-\delta(iT + p - s))) \sum_{k=0}^{i-1} \exp(-(\omega - \delta)Tk)$$

$$\leq \exp(-\omega iT) \varphi(p) + \gamma \exp(\delta T) \sup_{p \leq s \leq p+iT} (y(s) \exp(-\delta(iT + p - s))) \frac{1 - \exp(-(\omega - \delta)Ti)}{1 - \exp(-(\omega - \delta)T)}$$

Consequently, the following inequality holds for all $i \geq 0$ and $p \geq 0$:

$$\varphi(iT + p) \leq \exp(-\omega iT) \varphi(p) + \gamma \exp(\delta T) \sup_{p \leq s \leq p+iT} (y(s) \exp(-\delta(iT + p - s))) \frac{1 - \exp(-(\omega - \delta)Ti)}{1 - \exp(-(\omega - \delta)T)} \tag{A.6}$$

Inequality (4.94) implies that $\varphi(p) \leq M \exp(-\sigma p)\varphi(0) + \gamma \sup_{0 \leq s \leq p}(y(s))$, which combined with (A.6) gives for all $i \geq 0$ and $p \geq 0$:

$$\varphi(iT + p) \leq M \exp(-\omega iT - \sigma p) \varphi(0) + \gamma \exp(-\omega iT) \sup_{0 \leq s \leq p}(y(s))$$

$$+ \gamma \frac{\exp(\delta T)(1 - \exp(-(\omega - \delta)Ti))}{1 - \exp(-(\omega - \delta)T)} \sup_{p \leq s \leq iT + p} (y(s) \exp(-\delta(iT + p - s)))$$

$$\leq M \exp(-\omega iT - \sigma p) \varphi(0) + \gamma \exp(-\omega iT) \sup_{0 \leq s \leq p}(y(s) \exp(\delta s))$$

$$+ \gamma \frac{\exp(\delta T)(1 - \exp(-(\omega - \delta)Ti))}{1 - \exp(-(\omega - \delta)T)} \sup_{p \leq s \leq iT + p} (y(s) \exp(-\delta(iT + p - s)))$$

$$\leq M \exp(-\omega iT - \sigma p) \varphi(0) + \gamma \exp(-(\omega - \delta)iT) \exp(\delta p) \sup_{0 \leq s \leq p}(y(s) \exp(-\delta(iT + p - s)))$$

$$+ \gamma \frac{\exp(\delta T)(1 - \exp(-(\omega - \delta)Ti))}{1 - \exp(-(\omega - \delta)T)} \sup_{p \leq s \leq iT + p} (y(s) \exp(-\delta(iT + p - s)))$$

Using the above inequality in conjunction with the fact that $\omega \in (0,\sigma]$, we obtain for all $i \geq 0$ and $p \in [0,T]$:

$$\varphi(iT+p) \leq M \exp(-\omega(iT+p))\varphi(0) + \frac{\gamma \exp(\delta T)}{1-\exp(-(\omega-\delta)T)} \sup_{0 \leq s \leq iT+p}\left(y(s)\exp(-\delta(iT+p-s))\right) \tag{A.7}$$

Since (A.7) holds for all $i \geq 0$ and $p \in [0,T]$, we conclude that the following estimate holds for all $t \geq 0$:

$$\varphi(t) \leq M \exp(-\omega t)\varphi(0) + \frac{\gamma \exp(\delta T)}{1-\exp(-(\omega-\delta)T)} \sup_{0 \leq s \leq t}\left(y(s)\exp(-\delta(t-s))\right) \tag{A.8}$$

Definitions (A.3) guarantee that $\frac{\exp(\delta T)}{1-\exp(-(\omega-\delta)T)} = \frac{\mu}{1-\lambda\mu}$. Inequality (4.95) is a direct consequence of the previous equality, the fact that $\delta \in (0,\omega)$ and inequalities (A.2), (A.8).

The proof is complete. ◁

**Proof of Lemma 4.3:** Let $\varepsilon > 0$, $\lambda \in (0,1)$ (arbitrary) be given. Define

$$T = \sigma^{-1} \ln(M/\lambda), \quad \omega := \sigma - T^{-1}\ln(M), \quad \delta := \min\left(T^{-1}\ln(1+\varepsilon), \omega\right) \tag{A.9}$$

Notice that definitions (A.9) guarantee that $\lambda = M\exp(-\sigma T) = \exp(-\omega T) < 1$ and that $\omega \in (0,\sigma]$, $\delta \in (0,\omega]$. Using (4.102), we obtain for all $p \geq 0$:

$$q((i+1)T+p) \leq \max\left(\lambda q(iT+p), \gamma \sup_{p+iT \leq s \leq p+(i+1)T}(y(s))\right), \text{ for all } i \geq 0 \tag{A.10}$$

Applying (A.10) repeatedly, we get:

$$q(iT+p) \leq \max\left(\lambda^i q(p), \gamma \max_{k=0,\ldots i-1}\left(\lambda^{i-1-k} \sup_{p+kT \leq s \leq p+(k+1)T}(y(s))\right)\right), \text{ for all } i \geq 1 \text{ and } p \geq 0 \tag{A.11}$$

Using the facts that $\lambda = \exp(-\omega T)$ and $\delta \in (0,\omega]$, we get for all $i \geq 1$ and $p \geq 0$:

$$q(iT+p) \leq \max\left(\exp(-\omega iT)q(p), \gamma \exp(\omega T)\max_{k=0,\ldots i-1}\left(\exp(-\omega(i-k)T)\sup_{p+kT\leq s\leq p+(k+1)T}(y(s))\right)\right)$$

$$\leq \max\left(\exp(-\omega iT)q(p), \gamma \exp(\omega T)\max_{k=0,\ldots i-1}\left(\exp(-(\omega-\delta)(i-k)T)\sup_{p+kT\leq s\leq p+(k+1)T}(y(s)\exp(-\delta(iT+p-kT-p)))\right)\right)$$

$$\leq \max\left(\exp(-\omega iT)q(p), \gamma \exp(\delta T)\max_{k=0,\ldots i-1}\left(\exp(-(\omega-\delta)(i-1-k)T)\sup_{p+kT\leq s\leq p+(k+1)T}(y(s)\exp(-\delta(iT+p-s)))\right)\right)$$

$$\leq \max\left(\exp(-\omega iT)q(p), \gamma \exp(\delta T)\max_{k=0,\ldots i-1}\left(\sup_{p+kT\leq s\leq p+(k+1)T}(y(s)\exp(-\delta(iT+p-s)))\right)\right)$$

Consequently, the following inequality holds for all $i \geq 0$ and $p \geq 0$:

$$q(iT+p) \leq \max\left(\exp(-\omega iT)q(p), \gamma \exp(\delta T)\sup_{p \leq s \leq p+iT}(y(s)\exp(-\delta(iT+p-s)))\right) \tag{A.12}$$

Inequality (4.102) implies that $q(p) \leq \max\left(M\exp(-\sigma p)q(0), \gamma \sup_{0 \leq s \leq p}(y(s))\right)$, which combined with (A.12) and the facts that $\omega \in (0,\sigma]$, $\delta \in (0,\omega]$, gives for all $i \geq 0$ and $p \geq 0$:

$$q(iT+p)$$
$$\leq \max\left( M\exp(-\omega(iT+p))q(0), \gamma\exp(\delta p)\sup_{0\leq s\leq p}(y(s)\exp(-\delta(iT+p-s))), \gamma\exp(\delta T)\sup_{p\leq s\leq p+iT}(y(s)\exp(-\delta(iT+p-s)))\right)$$
(A.13)

Using (A.13), we obtain for all $i \geq 0$ and $p \in [0,T]$:

$$q(iT+p) \leq \max\left( M\exp(-\omega(iT+p))q(0), \gamma\exp(\delta T)\sup_{0\leq s\leq p+iT}(y(s)\exp(-\delta(iT+p-s)))\right) \quad (A.14)$$

Since (A.14) holds for all $i \geq 0$ and $p \in [0,T]$, we conclude that the following estimate holds for all $t \geq 0$:

$$q(t) \leq \max\left( M\exp(-\omega t)q(0), \gamma\exp(\delta T)\sup_{0\leq s\leq t}(y(s)\exp(-\delta(t-s)))\right) \quad (A.15)$$

Definitions (A.9) guarantee that $\exp(\delta T) \leq 1+\varepsilon$. Inequality (4.103) is a direct consequence of the previous inequality, the fact that $\delta \in (0, \omega]$ and inequality (A.15). The proof is complete. ◁